\documentclass [11pt] {article}
\usepackage[T1]{fontenc}
\usepackage{lmodern} 
\usepackage{textcomp}
\usepackage{amsmath}
\usepackage{amsfonts}
\usepackage{amsthm}
\usepackage{amssymb}
\usepackage{graphicx}
\usepackage{multicol}
\usepackage{color}
\usepackage{wrapfig}
\usepackage{subfig}
\usepackage[mathscr]{euscript}
\usepackage[all]{xy}
\usepackage{comment}
\usepackage{enumerate}
\usepackage{manfnt} 
\usepackage{array}
\usepackage{tikz-cd}
\usepackage[hidelinks]{hyperref}
\usepackage{cleveref}
\usepackage{aliascnt}
\usepackage{float}
\usepackage{caption}
\usepackage{microtype}
\usepackage{tabularx}
\usepackage{eso-pic}
\usetikzlibrary{fit, backgrounds}
\usetikzlibrary{calc}
\usetikzlibrary{decorations.markings} 
\usetikzlibrary{arrows.meta}           

\DeclareUnicodeCharacter{221E}{\ensuremath{\infty}}

\newtheoremstyle{nonitalic}
  {3pt}        
  {3pt}        
  {}           
  {}           
  {\bfseries}  
  {.}          
  {.5em}       
  {}           

\theoremstyle{plain}

\newtheorem{theorem}{Theorem}[section]
\crefname{theorem}{Theorem}{Theorems}
\Crefname{theorem}{Theorem}{Theorems}

\crefname{mainthm}{Theorem}{Theorems}
\Crefname{mainthm}{Theorem}{Theorems}

\newaliascnt{proposition}{theorem}
\newtheorem{proposition}[proposition]{Proposition}
\aliascntresetthe{proposition}
\crefname{proposition}{Proposition}{Propositions}
\Crefname{proposition}{Proposition}{Propositions}

\newaliascnt{corollary}{theorem}
\newtheorem{corollary}[corollary]{Corollary}
\aliascntresetthe{corollary}
\crefname{corollary}{Corollary}{Corollaries}
\Crefname{corollary}{Corollary}{Corollaries}

\newaliascnt{lemma}{theorem}
\newtheorem{lemma}[lemma]{Lemma}
\aliascntresetthe{lemma}
\crefname{lemma}{Lemma}{Lemmas}
\Crefname{lemma}{Lemma}{Lemmas}

\newaliascnt{conjecture}{theorem}
\newtheorem{conjecture}[conjecture]{Conjecture}
\aliascntresetthe{conjecture}
\crefname{conjecture}{Conjecture}{Conjectures}
\Crefname{conjecture}{Conjecture}{Conjectures}

\theoremstyle{nonitalic}

\newaliascnt{definition}{theorem}
\newtheorem{definition}[definition]{Definition}
\aliascntresetthe{definition}
\crefname{definition}{Definition}{Definitions}
\Crefname{definition}{Definition}{Definitions}

\newaliascnt{example}{theorem}

\aliascntresetthe{example}
\crefname{example}{Example}{Examples}
\Crefname{example}{Example}{Examples}

\newaliascnt{rem}{theorem}
\newtheorem{rem}[rem]{Remark}
\aliascntresetthe{rem}
\crefname{rem}{Remark}{Remarks}
\Crefname{rem}{Remark}{Remarks}

\newtheorem*{rem*}{Remark}

\newcommand{\R}{\mathbb{R}}

\newcommand{\cC}{\mathcal{C}}

\newcommand{\cS}{\mathcal{S}}
\newcommand{\cE}{\mathcal{E}}

\newcommand{\cU}{\mathcal{U}}

\newcommand{\lf}{\left}
\newcommand{\ri}{\right}
\newcommand{\norm}[1]{\left\vert{#1}\right\vert} 

\newcommand{\Hom}{\operatorname{Hom}}
\newcommand{\Psh}{\operatorname{PSh}}
\newcommand{\Shv}{\operatorname{Shv}}

\newcommand{\Top}{\operatorname{Top}}
\newcommand{\Strat}{\operatorname{CStrat}}

\newcommand{\Fun}{\operatorname{Fun}}

\newcommand{\Set}{\operatorname{Set}}
\newcommand{\sSet}{\operatorname{sSet}}
\newcommand{\Sing}{\operatorname{Sing}}

\newcommand{\Sd}{\operatorname{Sd}}

\newcommand{\Cat}{\operatorname{Cat}_\infty}
\newcommand{\cst}{\operatorname{cst}}
\newcommand{\colim}[2]{\mathop{\operatorname{colim}}\limits_{#1} #2}
\newcommand{\limi}[2]{\mathop{\operatorname{lim}}\limits_{#1} #2}

\newcommand{\exit}[2]{\operatorname{Exit}^{#1}(#2)}

\newcommand{\pos}[1]{\operatorname{B}({#1})}

\newcommand{\Pre}{\operatorname{PreOrd}}
\newcommand{\Arr}{\operatorname{Arr}}

\tikzset{
  midarrow/.style={
    postaction={decorate},
    decoration={
      markings,
      mark=at position 0.5 with {\arrow{>}}
    }
  }
}

\makeatother

\begin{document}

\title{A stratified version of geometric realizations}
\date{2026}

\author{Rodrigue Haya Enriquez}


\maketitle

\begin{abstract}
Given a simplicial set $X$, let $\pos{X}$ denote its poset of nondegenerate simplices. We define a natural stratification $\norm{X}\rightarrow\pos{X}$. When $X$ is locally finite, we prove that this stratification is conical. Moreover, if $X$ is finite, we prove that the exit path category $\exit{\pos{X}}{\norm{X}}$ satisfies the geometric realization formula $\colim{\sigma:\Delta[p]\rightarrow X}{\exit{\pos{p}}{\Delta^p}}\simeq\exit{\pos{X}}{\norm{X}}$ where the colimit is indexed by the category of simplices of $X$.
\end{abstract}

\begingroup
\small
\setlength{\parskip}{2pt}
\linespread{0.9}\selectfont
\tableofcontents
\endgroup


\section{Introduction}

Let $T$ be a topological space and $P$ be a poset. Informally, $P$-stratification on $T$ is a decomposition $T=\bigsqcup\limits_{u\in P}T_u$ into \emph{strata}, with mild compatibility conditions. Following \cite[Appendix A]{ha}, this is defined as a continuous map $T\rightarrow P$, where $P$ is a poset endowed with its Alexandrov topology. 

Now let $X$ be a simplicial set and $\pos{X}$ be the poset whose elements are the nondegenerate simplices ordered by face relations. The geometric realization $\norm{X}$ of $X$ admits a stratification $\pi_X:\norm{X}\rightarrow\pos{X}$ sending a point $x$ to the unique nondegenerate simplex whose interior contains $x$. 

The present paper is mainly interested in the \emph{exit path category} of this stratification. Provided it exists, it is an $\infty$-category whose objects are the points of $\norm{X}$ and the morphisms are paths in $\norm{X}$ than can go from lower dimensional simplices to upper dimensional simplices but not the other way around.

In \cite[Appendix A.6]{ha}, Lurie introduces the \emph{conicality} condition on a stratification and then shows that conically stratified spaces admit an exit path category (\cite[Theorem A.6.4]{ha}). Our first result is that the stratification $\pi_X:\norm{X}\rightarrow\pos{X}$ is conical.

\begin{theorem}\label{main_strat_A}
Let $X$ be a locally finite simplicial set. Then the stratification $\pi_X:\norm{X}\rightarrow\pos{X}$ is conical.
\end{theorem}

Assuming that $X$ is locally finite, \cite[Theorem A.6.4]{ha} and \Cref{main_strat_A} ensure that the stratification $\pi_X:\norm{X}\rightarrow\pos{X}$ admits an exit path category denoted $\exit{\pos{X}}{\norm{X}}$. We mention that \Cref{main_strat_A} was proved by Lurie for simplicial complexes (\cite[Proposition A.6.8]{ha}). 

Next, for $p\geq 0$, the geometric $p$-simplex $\Delta^p$ admits a stratification by the poset $\pos{p}$ of its faces. The assignment $p\mapsto\pos{p}$ extends to a cosimplicial space $\pos{-}:\Delta\rightarrow\Top$. Let $\Delta(X)$ be the category of simplices of $X$. Similarly to the geometric realization formula

\[
\colim{(\sigma:\Delta[p]\rightarrow X)\in\Delta(X)}{\Delta^p}\cong\norm{X},
\]
we have a homeomorphism (\Cref{colim_B(X)})
\[
\colim{(\sigma:\Delta[p]\rightarrow X)\in\Delta(X)}{\pos{p}}\cong\pos{X}.
\]

This raises the question whether such a colimit passes to exit path categories. This is not trivial, since the exit path functor does not preserve arbitrary colimits and degeneracies can create nontrivial mapping spaces. Nonetheless, we provide a positive answer to this question when $X$ is finite.

\begin{theorem}\label{main_strat_B}
Let $X$ be a finite simplicial set. Then there is an equivalence of $\infty$-categories
\[
\colim{(\sigma:\Delta[p]\rightarrow X)\in\Delta(X)}{\exit{\pos{p}}{\Delta^p)}}\simeq\exit{\pos{X}}{\norm{X}}.
\]
\end{theorem}

\Cref{main_strat_B} gives a tool to compute (co)limits indexed by exit path categories. Suppose that $S\rightarrow P$ is a conical stratification and $X$ a finite simplicial set providing a a refinement 
\[
\begin{tikzcd}
\norm{X}\arrow[d,"\pi_X"']\arrow[r,"\cong"] & S\arrow[d,"\pi"] \\
\pos{X}\arrow[r] & P
\end{tikzcd}.
\]
For example, $S$ could be a manifold and $X$ is obtained from a triangulation. Then \Cref{main_strat_A} and \Cref{main_strat_B} can be used to provide an explicit ordinary category $\cE$ together with a localization functor $\cE\rightarrow \exit{P}{S}$. We refer to \cite[Section 6-7]{main_haya} for a detailed application of this idea.

In the present paper, we start by reviewing stratified spaces, exit path categories, constructible sheaves and the exodromy equivalence based on \cite[Appendix A]{ha} and \cite{jansen}. Since we do not need the results from these references in their full generality, and in order to keep this document not too involved, we will restrict our attention to locally finite CW-complexes and space-valued sheaves on them. As we will be interested in stratifications on geometric realizations of finite simplicial sets, this simplification causes no issue. In the third section, we define for a simplicial set $X$ the poset $\pos{X}$ and the stratification $\pi_X:\norm{X}\rightarrow\pos{X}$ sketched above. In the fourth section, we prove that for a locally finite simplicial set $X$, this stratification is conical. The proof is unfortunately quite involved and technical, which accounts for its length. Combined with \cite[Theorem A.6.4]{ha}, this result allows us to define explicitly the exit path category relative to this stratification. Based on this and the material recalled in the second section, we finally prove \Cref{main_strat_B}.

\subsection*{Conventions}

Throughout the text, $\infty$-categories are modelled by \emph{quasicategories}. We denote by $\cS$ the $\infty$-category of spaces and $\Cat$ the $\infty$-category of (small) $\infty$-categories.

For $p\geq 0$, we write $\Delta[p]$ for the representable simplicial set of dimension $p$ and $\Delta^p$ the geometric $p$-simplex. Given a simplicial set $X$ and a simplex $\sigma:\Delta[p]\rightarrow X$, we write $\Delta[\sigma]:=\Delta[p]$ and $\Delta^\sigma:=\Delta^p$.  We also write $\Delta(X)$ for its category of simplices and $\Delta^{inj}(X)$ for the subcategory having the same objects but only face maps as morphisms.

We write $\Arr(\Top)$ for the arrow category of $\Top$, that is, the category whose objects are continuous maps.

For a category $\cC$, we write $\cC^\lhd$ (resp. $\cC^\rhd$) the category obtained by adjoining an initial object $-\infty$ (resp. a terminal object $\infty$).

\subsection*{Acknowledgements}
This work was realized during the course of the author's PhD program.

He thanks his advisors Pascal Lambrechts and Nathalie Wahl for their careful reading and for providing valuable comments on earlier versions of this document.

He also thanks Colin Fourel for answering questions about exit path categories during the
"Young Topologists Meeting 2026" and suggesting the paper \cite{jansen}.

The author was supported by the "Fonds de la Recherche Scientifique–FNRS" under the Aspirant Fellowship FC 49827. 

\subsection*{Use of generative AI}
During the preparation of this document, generative AI has been used for skimming more rapidly the relevant literature, reviewing early versions by pointing out mistakes and inconsistencies and assisting with the document layout. The models used were Gemini 3.5 and ChatGPT-5.6.

\section{Exit path categories and constructible sheaves}
Let $T$ be a topological space. We are interested in \emph{stratifications} on $T$, that is a decomposition $T=\bigsqcup\limits_{u\in P}T_u$ into \emph{strata}, with mild compatibility conditions. One way of formalizing this idea is done in \cite[Appendix A.5]{ha}. Namely, a stratification on $T$ is simply given by a map $\pi:T\rightarrow P$ where $P$ is a poset and the strata are the preimages $\{T_u:=\pi^{-1}(u)\}_{u\in P}$. 

This section is a brief review of \cite[Appendix A]{ha}. In addition, we also recall some results from \cite{jansen}, extracting the relevant material for the present paper.

\begin{definition}
Let $P$ be a poset. The \emph{Alexandroff topology} on $P$ is the topology where a subset $U\subset P$ is open if for any $u\in U$ and $v\in P$ such that $u\leq v$, we have $v\in U$.
\end{definition}

In particular, for any $u\in P$, the subsets $P_{u\leq}:=\{v\in P,u\leq v\},P_{u<}:=\{v\in P,u<v\}$ are open. 

\begin{definition}
Given a topological space $S$ and a poset $P$, a \emph{$P$-stratification on $S$} consists of a continuous map $\pi:S\rightarrow P$. A \emph{stratified space} is a topological space endowed with a $P$-stratification for some poset $P$.

Given stratified spaces $(\pi:S\rightarrow P),(\psi:T\rightarrow Q)$, a \emph{stratified map} consists of continuous maps $f^\sharp:S\rightarrow T,f^\flat:P\rightarrow Q$ making the following square commute:
\[
\begin{tikzcd}
S\arrow[d,"\pi"']\arrow[r,"f^\sharp"] & T\arrow[d,"\psi"] \\
P\arrow[r,"f^\flat"'] & Q
\end{tikzcd}.
\]

\end{definition}

Given a stratified space $\pi:S\rightarrow P$ satisfying some conditions, we will be interested in its \emph{exit path category}. Informally, this consists of an $\infty$-category whose objects are the points of $S$ and the morphisms are paths only allowed to go from lower to upper strata.

\begin{definition}
Given a poset $P$, let $P^{\lhd}$ be the poset obtained by adjoining a minimal element $-\infty$ to $P$, that is $P^{\lhd}=\{-\infty\}\sqcup P$ and $-\infty<u$ for any $u\in P$.
\end{definition}

\begin{definition}\cite[Definition A.5.3]{ha}
Given a topological space $S$, the open cone $C(S)$ is the topological space whose underlying set is $\{\ast\}\sqcup S\times \R_{>0}$. The \emph{teardrop topology} on $C(S)$ is the topology  for which $U\subset C(S)$ is open if and only if $U\cap (S\times\R_{>0})$ is open, and whenever $\ast\in U$, there exists $\epsilon>0$ such that $S\times ]0,\epsilon[\subset U$.
\end{definition}

\begin{lemma}\label{teardrop_topology}
Let $S$ be a compact topological space. Then the cone $C(S)$ endowed with the teardrop topology is homeomorphic to the pushout $S\times\R_{\geq 0}\sqcup_{S\times\{0\}}\ast$. 
\end{lemma}

\begin{proof}
We prove that a subset $U\subset C(S)$ is open with respect to the teardrop topology if and only if it is open with respect to the final topology on $C(S)$ defined by the quotient map $q:S\times\R_{\geq 0}\rightarrow C(S)$,
\[
q(x,t)=
\lf\{
\begin{array}{ll}
\ast & t=0, \\
(x,t)\in S\times\R_{>0} & t>0.
\end{array}
\ri.
\] 

Assume that $U$ is open with respect to the teardrop topology. If $\ast\not\in U$, then $U=U\cap(S\times\R_{>0})$ is open in $S\times\R_{>0}$. Thus $q^{-1}(U)=U\cap (S\times\R_{>0})$ is open in $S\times\R_{\geq 0}$. By definition, this implies that $U$ is open with respect to the final topology. If $\ast\in U$, there exists $\epsilon>0$ such that $S\times ]0,\epsilon[\subset U$. Decomposing $U$ as 
\[
U=\lf(\{\ast\}\sqcup (S\times ]0,\epsilon[)\ri)\cup \lf(U\cap S\times\R_{>0}\ri),
\]
we get 
\begin{align*}
q^{-1}(U) &= q^{-1}\lf(\{\ast\}\sqcup (S\times ]0,\epsilon[)\ri)\cup q^{-1}\lf(U\cap S\times\R_{>0}\ri) \\
& =(S\times [0,\epsilon[)\cup q^{-1}\lf(U\cap S\times\R_{>0}\ri)
\end{align*}
which is open in $S\times\R_{\geq 0}$. Again, this implies that $U$ is open with respect to the final topology defining the pushout.

Conversely, assume that $U$ is open with respect to the final topology. If $\ast\not\in U$, then $q^{-1}(U)=U\cap (S\times\R_{>0})$ is open in $S\times\R_{\geq 0}$ and thus is open in $S\times\R_{>0}$. Otherwise, $S\times\{0\}\subset q^{-1}(U)$. Then there exists a family of open subsets $\{V_\alpha\times [0,\epsilon_\alpha[,V_\alpha\subset S,\epsilon_\alpha>0\}_{\alpha\in I}$ where the family $\{V_\alpha\}_{\alpha\in I}$ forms an open cover of $S$ and $V_\alpha\times [0,\epsilon_\alpha[\subset q^{-1}(U)$ for all $\alpha\in I$. Since $S$ is compact, there exists a finite subset $\{V_{\alpha_i}\}_{i=1,\cdots,N}\subset\{V_\alpha\}_{\alpha\in I}$ providing an open cover of $S$. Defining $\epsilon:=\min\limits_{i=1,\cdots,N}\epsilon_{\alpha_i}$, the subset $S\times [0,\epsilon[$ is open in $S\times\R_{\geq 0}$ and 
\[
S\times [0,\epsilon[\subset\bigcup_{i=1,\cdots,N}(V_{\alpha_i}\times [0,\epsilon_{\alpha_i}[)\subset q^{-1}(U).
\] 
From this, we deduce that $S\times ]0,\epsilon[\subset U$. In addition, $U\cap(S\times\R_{>0})=q^{-1}(U\cap (S\times R_{>0}))$ is open in $S\times\R_{\geq 0}$ hence $U\cap(S\times\R_{>0})$ is open in $S\times\R_{>0}$. Combining both observations yields that $U$ is open with respect to the teardrop topology.

\end{proof}

\begin{definition}
Let $\pi:S\rightarrow P$ be a stratified space. We define a $P^{\lhd}$-stratification on the cone $C(S)$ by

\[
x\mapsto
\left\{
\begin{array}{ll}
-\infty & \text{if }x=\ast. \\
\pi(s) & x=(s,t)\in S\times\R_{>0}.
\end{array}
\right.
\]
\end{definition}

For the following definition, first observe that given a topological space $Y$ and a stratified space $S\rightarrow P$, the product $Y\times S$ is also $P$-stratified by the composition $Y\times S\xrightarrow{\pi_S} S\rightarrow P$. Moreover for any $a\in P$, there is an isomorphism of posets $(P_{a<})^{\lhd}\cong P_{a\leq}$.

\begin{definition}
Let $\pi:S\rightarrow P$ be a stratified space. We say that this stratification is \emph{conical} if for any point $x\in S$ lying over an element $a\in P$, there exists a $P_{a<}$-stratified space $L$, a topological space $U$ and a commutative square

\[
\begin{tikzcd}
U\times C(L)\arrow[r,hook,"\iota"]\arrow[d] & S\arrow[d] \\
P_{a\leq}\arrow[r,hook] & P 
\end{tikzcd}
\]
where the map $\iota:U\times C(L)\hookrightarrow S$ is an open embedding with $x=\iota(u,\ast)$ for some $u\in U$.
\end{definition}

\begin{definition}
The category $\Strat$ has objects conically stratified spaces and morphisms stratified maps between them. 
\end{definition}

Let $\Delta^p=\{(t_0,\cdots,t_p)\in [0,1]^{p+1},\sum_{i=0}^p t_i=1\}$ be the geometric $p$-simplex and $[p]=\{0\leq 1\leq\cdots\leq p\}$. There is a conical stratification $\lambda^p:\Delta^p\rightarrow [p]$ sending a point $t=(t_0,\cdots,t_p)$ to $\max\{k\in [p],t_k\neq 0\}$. This defines a functor $(\lambda^\bullet:\Delta^\bullet\rightarrow [\bullet]):\Delta\rightarrow\Strat$. 

For a poset $P$, its nerve $N(P)\in\sSet$ is given in degree $p$ by the set of morphisms of posets $\{f:[p]\rightarrow P\}$. Such a morphism $f$ defines a continuous maps with respect to the Alexandrov topologies and precomposition with the map $\lambda^p:\Delta^p\rightarrow [p]$ yields a continuous map $f\circ\lambda^p:\Delta^p\rightarrow P$, that is a $p$-simplex of the singular set of $P$. By naturality of the map $\lambda^p$, this defines a morphism $N(P)\rightarrow\Sing(P)$.

Following Lurie, we have the following model for the \emph{exit path category} of a stratified space.
\begin{proposition}\cite[Theorem A.6.4]{ha}
Let $\pi:S\rightarrow P$ be a conically stratified space. Then the pullback of simplicial sets
\[
\begin{tikzcd}
N(P)\times_{\Sing(P)}\Sing(S)\arrow[d]\arrow[r]\arrow[dr, phantom, "\lrcorner", very near start] & \Sing(S)\arrow[d] \\
N(P)\arrow[r] & \Sing(P)
\end{tikzcd}
\]
is an $\infty$-category. 
\end{proposition}
If $\pi:S\rightarrow P$ is a conically stratified space, we will write $\exit{P}{S}:=N(P)\times_{\Sing(P)}\Sing(S)$. 

\begin{rem}
Unwinding the definition, we get that $\exit{P}{S}_p\cong\Hom_{\Strat}((\Delta^p,[p]),(S,P))$. In other words, any $p$-simplex is given by some stratified map
\[
\begin{tikzcd}
\Delta^p\arrow[d]\arrow[r] & S\arrow[d] \\
{[p]}\arrow[r] & P
\end{tikzcd}.
\]
\end{rem}

The above construction extends to a functor $\exit{-}{-}:\Strat\rightarrow\Cat$. The next step is to relate such exit path categories to \emph{constructible sheaves} on $S$.

Given a topological space $S$, we write $\Shv(S)$ for the $\infty$-category of space-valued sheaves on $S$.

\begin{definition}
Let $S$ be a topological space and $F\in\Shv(S)$.
\begin{enumerate}
\item $F$ is constant if it is in the essential image of the pullback functor $p^*:\cS\simeq\Shv(\ast)\rightarrow\Shv(S)$ induced by the map $p:S\rightarrow\ast$.
\item $F$ is locally constant if there exists an open cover $\cU$ of $S$ such that for every $U\in\cU$, the restriction $j_U^*(F)\in\Shv(U)$ induced by the inclusion $j_U:U\hookrightarrow S$ is a constant sheaf.
\end{enumerate}
We write $\Shv^{loc}(S)\subset\Shv(S)$ for the full subcategory of locally constant sheaves.
\end{definition}

For simplicity and since it is enough for our purposes, we will assume that the prescribed topological space $S$ is a locally finite CW-complex but the result holds in more generality, see \cite[Appendix A]{ha}. Given such a topological space and a locally constant sheaf $F$, we get a well-defined functor $F_{-}:\Sing(S)\rightarrow\cS$ sending a point $x$ to its stalk $F_x$. By straightening-unstraightening, this defines a Kan fibration $\cE\rightarrow \Sing(S)$. Conversely, given such a fibration, we get a locally constant sheaf 
\[
\psi_S(\cE):\text{Open}(S)^{op}\rightarrow\cS,U\mapsto\cS_{/\Sing(S)}(\Sing(U),\cE).
\] 

\begin{theorem}\cite[Theorem A.4.19]{ha}
Let $S$ be a locally finite CW-complex. Then there is an equivalence of $\infty$-categories 
\[
\psi_S:\cS_{/\Sing(S)}\xrightarrow{\simeq}\Shv^{loc}(S).
\]
\end{theorem}

The next step is to state a variant of this result for stratified spaces.

\begin{definition}
Let $\pi:S\rightarrow P$ be a stratified space. A sheaf $F\in\Shv(S)$ is \emph{constructible} if for any $p\in P$, the inclusion $j_p:S_p\hookrightarrow S$ induces a locally constant sheaf $j_p^*(F)\in\Shv(S_p)$. We write $\Shv^P(S)\subset\Shv(S)$ for the full subcategory of constructible sheaves.
\end{definition}

Now assume that the stratification $\pi:S\rightarrow P$ is conical. By \cite[Theorem A.6.4]{ha}, the simplicial set $\exit{P}{S}$ is an $\infty$-category. Consider now the $\infty$-category $\text{LFib}_{/\exit{P}{S}}$ of left fibrations over $\exit{P}{S}$. Following \cite[Construction A.9.2]{ha}, the functor 
\[
\psi_S^P:\text{LFib}_{/\exit{P}{S}}\rightarrow\Psh(\text{Open}(S))
\] 
is defined by 
\[
\psi_S^P(\cE):=\lf(U\mapsto \Fun_{/\exit{P}{S}}(\exit{P}{U},\cE)\ri).
\] 
Since $\cE\rightarrow\exit{P}{S}$ is a left fibration, the mapping functor is indeed a Kan complex. We can now state Lurie's version of the exodromy equivalence, at least in a particular case.

\begin{theorem}\cite[Theorem A.9.3]{ha}
Let $S$ be a locally finite CW-complex, $P$ a finite poset and $\pi:S\rightarrow P$ a conical stratification. Then there is an equivalence of $\infty$-categories 
\[
\psi^P_S:\text{LFib}_{/\exit{P}{S}}\xrightarrow{\simeq}\text{Shv}^P(S).
\]
\end{theorem}

In particular, composing with the straightening-unstraightening equivalence $\text{LFib}_{/\exit{P}{S}}\simeq\Fun(\exit{P}{S},\cS)$ provides an equivalence of $\infty$-categories $\Fun(\exit{P}{S},\cS)\simeq\Shv^P(S)$. The upshot is that constructible sheaves can be recovered from the exit path category $\exit{P}{S}$ under the previous assumptions. Nonetheless, we will instead be interested in proving properties about exit path categories by exploiting $\infty$-categories of constructible sheaves and the exodromy equivalence. 

This leads to results from \cite{jansen}. Before pursuing, we note that while Jansen's definition of exit path category differs from Lurie's definition introduced before, both definitions agree for conical stratifications on locally finite CW-complexes over a finite poset, see \cite[Example 3.13]{jansen}. Moreover, based on \cite[Theorem 3.2.3]{exodromy}, \cite[Remark 3.15]{jansen} points  out that given a stratified map $f:(S,P)\rightarrow (T,Q)$ where both $(S,P),(T,Q)$ are conical stratifications on locally finite CW-complexes over a finite poset, the induced functor $\exit{P}{S}\rightarrow\exit{Q}{T}$ respects exit path categories in the sense of \cite[Definition 3.12]{jansen}.

\begin{proposition}\label{jansen_2_20}\cite[Corollary 2.20]{jansen}
Consider a pullback square of topological spaces where each object is a locally finite CW-complex:
\[
\begin{tikzcd}
S'\arrow[d,"\iota'"']\arrow[r]\arrow[dr,"\lrcorner",phantom,very near start] & S\arrow[d,hook,"\iota"] \\
T'\arrow[r,"q"'] & T
\end{tikzcd}.
\]
Assume that 
\begin{enumerate}
\item $\iota:S\hookrightarrow T$ is a closed embedding.
\item $q:T'\rightarrow T$ is a proper map.
\item $q$ induces a homeomorphism $T'\setminus\iota'(S')\xrightarrow{\cong}T\setminus\iota(S)$.
\end{enumerate}
Then the induced square
\[
\begin{tikzcd}
\Shv(T)\arrow[d,"q^*"']\arrow[r,"\iota^*"] & \Shv(S)\arrow[d] \\
\Shv(T')\arrow[r,"(\iota')^*"'] & \Shv(S')
\end{tikzcd}
\]
is a pullback square in $\Cat$.
\end{proposition}

The following proposition is a reformulation of \cite[Proposition 3.16]{jansen} in a particular case. It allows us to compute a colimit of exit path categories based on the limit of constructible sheaves they classify.

\begin{proposition}\label{jansen_3_16}\cite[Proposition 3.16]{jansen}
Let $I$ be a small ordinary category and 
\[
F:I^\rhd\rightarrow\Strat, F(u)=\pi_u:S_u\rightarrow P_u
\]
a cocone where each $\pi_u:S_u\rightarrow P_u$ is a conical stratification of a locally finite CW-complex over a finite poset and $F(\infty):\pi_\infty:S_\infty\rightarrow P_\infty$ is the image of the terminal object $\infty\in I^\rhd$. Assume moreover that the cocone induces equivalences of $\infty$-categories
\begin{enumerate}
\item $\Shv(S_\infty)\xrightarrow{\simeq}\limi{u\in I}{\Shv(S_u)}$.
\item $\Shv^{P_\infty}(S_\infty)\xrightarrow{\simeq}\limi{u\in I}{\Shv^{P_u}(S_u)}$.
\end{enumerate}  
Then the induced functor $\colim{u\in I}{\exit{P_u}{S_u}}\rightarrow\exit{P_\infty}{S_\infty}$ is an equivalence of $\infty$-categories.
\end{proposition}

\section{The stratification $\pi_X:\norm{X}\rightarrow\pos{X}$}
In this section, for each simplicial set $X$, we investigate its poset of \emph{nondegenerate simplices}.

We start with the base case of representable simplicial sets. For each $p\geq 0$, consider the linearly ordered set $[p]=\{0\leq 1\leq\cdots\leq p\}$ and $\pos{p}$ the poset of nonempty subsets of $[p]$ ordered by inclusions. This defines an abstract simplicial complex $([p],\pos{p})$. The geometric realization $\norm{\Delta^{([p],\pos{p})}}$ of this abstract simplicial complex is identified with the topological space of maps $t:[p]\rightarrow [0,1]$ such that $\sum_{i\in [p]}t_i=1$ and $\text{Supp}(t)=\{i\in [p],t_i\neq 0\}\in\pos{p}$. This is nothing but the geometric simplex $\Delta^p=\{(t_0,\cdots,t_p)\in [0,1]^{p+1},\sum_{i=0}^p t_i=1\}$.

Following \cite[Definition A.6.7]{ha}, we get a stratification $\pi_p:\Delta^p\rightarrow \pos{p},t=(t_0,\cdots,t_p)\mapsto\text{Supp}(t)$. By \cite[Proposition A.6.8]{ha}, this stratification is conical.

\begin{lemma}
The stratification $\pi_p:\Delta^p\rightarrow\pos{p}$ extends to a functor $\pi_\bullet:\Delta\rightarrow\Strat$.
\end{lemma}

\begin{proof}
Any linearly ordered map $f:[p]\rightarrow [q]$ induces a morphism $\pos{f}:\pos{p}\rightarrow\pos{q}$ sending $S\subset [p]$ to its image $f(S)\subset [q]$. We check that such a morphism is compatible with the map $\Delta^f:\Delta^p\rightarrow\Delta^q$. It is enough to check this compatibility for cofaces and codegeneracies.
Consider a coface $d^i:[p-1]\rightarrow [p]$ for some $0\leq i\leq p$. Given $t=(t_0,\cdots,t_{p-1})\in\Delta^{p-1}$, we have $\Delta^f(t)=(t'_0,t'_1,\cdots,t'_{p})$ with $t'_j=t_j$ if $j<i$, $t'_i=0$ and $t'_j=t_{j-1}$ for $j>i$. Hence $t'_j\neq 0$ if and only if $j\in f(\text{Supp}(t))$.

Now consider a codegeneracy $s^i:[p+1]\rightarrow [p]$ for some $0\leq i\leq p$. Given $t=(t_0,\cdots,t_{p+1})\in\Delta^{p+1}$, we have $\Delta^f(t)=(t'_0,t'_1,\cdots,t'_{p})$ with $t'_j=t_j$ if $j<i$, $t'_i=t_i+t_{i+1}$ and $t'_j=t_{j+1}$ for $j>i$. For $j\neq i,i+1$, we thus have $j\in\text{Supp}(t)$ if and only if $s^i(j)\in\text{Supp}(t')$. Moreover, $i\in\text{Supp}(t')$ if and only if one of the elements $i,i+1$ is in $\text{Supp}(t)$. Combining both cases, we get that $\text{Supp}(t')=f(\text{Supp}(t))$.
\end{proof}

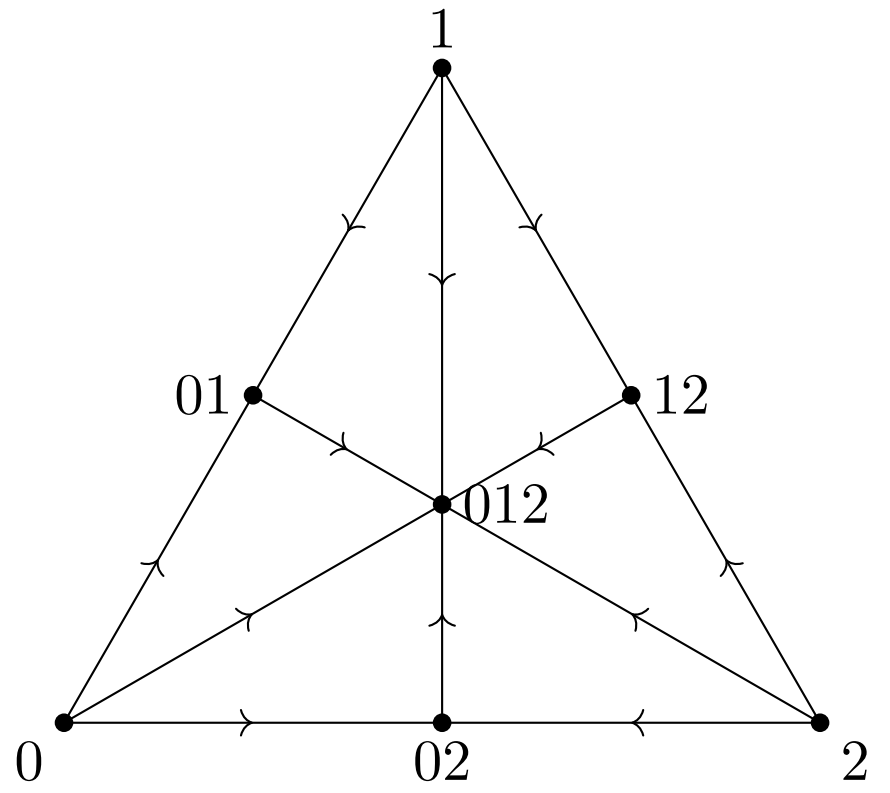
\begin{figure}[ht]
\centering
\begin{tikzpicture}[scale=2.5, line cap=round, line join=round]
  \tikzset{
    midarrow/.style={
      postaction={decorate},
      decoration={
        markings,
        mark=at position 0.5 with {\arrow{>}}
      }
    }
  }

  \coordinate (0) at (0,0);
  \coordinate (2) at (2,0);
  \coordinate (1) at (1,{sqrt(3)});

  \coordinate (01) at ($(0)!0.5!(1)$);
  \coordinate (02) at ($(0)!0.5!(2)$);
  \coordinate (12) at ($(1)!0.5!(2)$);

  \coordinate (012) at ($(0)!1/2!(1)!1/3!(2)$);

  \draw[midarrow] (0) -- (01);
  \draw[midarrow] (0) -- (02);
  \draw[midarrow] (1) -- (01);
  \draw[midarrow] (1) -- (12);
  \draw[midarrow] (2) -- (02);
  \draw[midarrow] (2) -- (12);

  \draw[midarrow] (01) -- (012);
  \draw[midarrow] (02) -- (012);
  \draw[midarrow] (12) -- (012);
  
  \draw[midarrow] (0) -- (012);
  \draw[midarrow] (1) -- (012);
  \draw[midarrow] (2) -- (012);

  \fill (0) circle (0.7pt);
  \fill (1) circle (0.7pt);
  \fill (2) circle (0.7pt);
  \fill (01) circle (0.7pt);
  \fill (02) circle (0.7pt);
  \fill (12) circle (0.7pt);
  \fill (012) circle (0.7pt);

  \node[below left]  at (0)   {$0$};
  \node[below right] at (2)   {$2$};
  \node[above]       at (1)   {$1$};

  \node[below]       at (02)  {$02$};
  \node[left]        at (01)  {$01$};
  \node[right]       at (12)  {$12$};

  \node[right] at (012) {$012$};
\end{tikzpicture}
\caption{Depicting the elements of $\pos{2}$ as the barycenters of $\Delta^2$.}
\end{figure}

Before generalizing this construction for arbitrary simplicial sets, we introduce a different point of view by considering the poset $\pos{p}$ as the poset of nondegenerate simplices of the standard simplicial set $\Delta[p]$ ordered by face relations. Namely, nondegenerate simplices of this simplicial set are cofaces $f:[k]\hookrightarrow [p]$ and such a coface uniquely determines a subset $S=f([k])\subset [p]$. The partial order is then given by $(f:[k]\hookrightarrow [p])\leq (g:[l]\hookrightarrow [p])$ whenever there exists a coface $\theta:[k]\rightarrow [l]$ such that $g\theta = f$.

This suggests a generalization for arbitrary simplicial sets.  

\begin{definition}
Let $X$ be a simplicial set. The poset $\pos{X}$ has as elements the nondegenerate simplices of $X$, with the partial order given by $\sigma\leq\tau$ whenever there exists a face map $f$ such that $\sigma = f(\tau)$.
\end{definition}

Given a simplex $\sigma:\Delta[p]\rightarrow X$, we will write $\pos{\sigma}:=\pos{p}$. 

\begin{lemma}\label{colim_preord}
There is a colimit of preordered sets $\colim{\sigma\in\Delta(X)}{\pos{\sigma}}\cong\pos{X}$.
\end{lemma}

\begin{proof}
Recall that for each simplex $\sigma$ of $X$, the Eilenberg-Zilberg lemma guarantees the existence of a unique nondegenerate simplex $ez(\sigma)$ for which $\sigma$ is a degeneracy.

For each simplex $\sigma:\Delta[p]\rightarrow X$, we define a morphism $\lambda_{\sigma}:\pos{\sigma}\rightarrow \pos{X}$ by 
\[
(f:[k]\hookrightarrow [p])\mapsto ez(\Delta[k]\xrightarrow{f} \Delta[p]\xrightarrow{\sigma} X).
\]

Given a morphism $\theta:[q]\rightarrow [p]$, we first check that $\lambda_{\sigma}\circ\theta = \lambda_{\sigma\theta}$. Let $(f:[k]\hookrightarrow [q])\in \pos{q}$. We have the commutative diagram

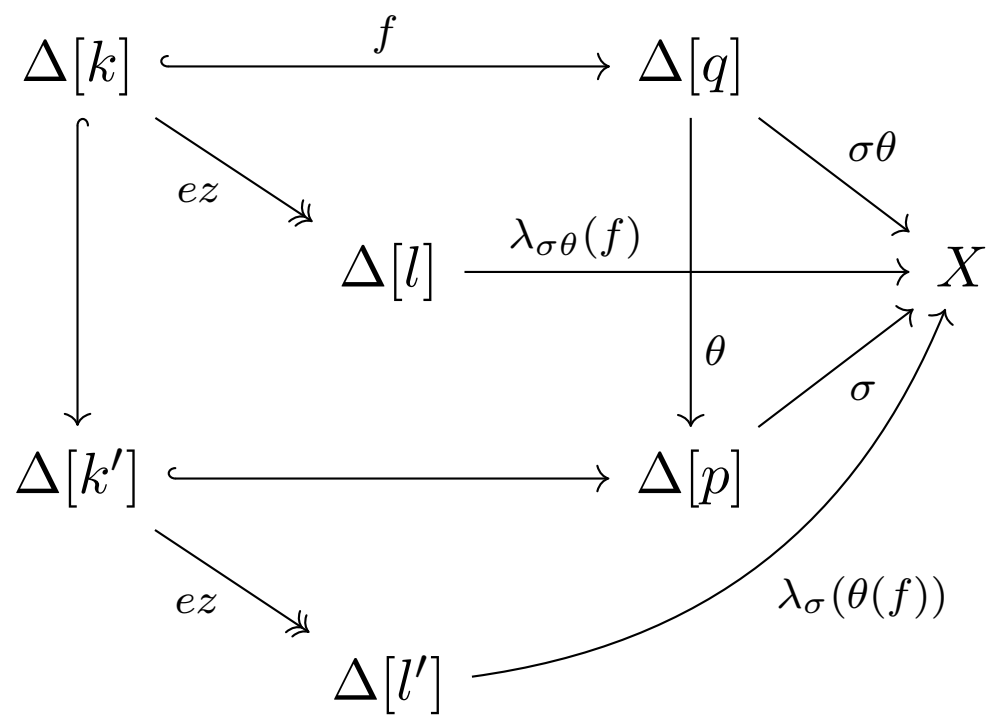
\[
\begin{tikzcd}
\Delta[k]\arrow[dr,two heads,"ez"']\arrow[rr,hook,"f"]\arrow[dd,hook] & & \Delta[q]\arrow[dd,"\theta",near end]\arrow[dr,"\sigma\theta"] & \\
& \Delta[l]\arrow[rr,"\lambda_{\sigma\theta}(f)",near start] & & X \\
\Delta[k']\arrow[rr,hook]\arrow[dr,two heads,"ez"'] & & \Delta[p]\arrow[ur,"\sigma"'] & \\
& \Delta[l']\arrow[uurr,bend right,"\lambda_{\sigma}(\theta(f))"'] & &
\end{tikzcd}.
\] 
We observe that the simplex $\Delta[k]\xrightarrow{f}\Delta[q]\xrightarrow{\sigma\theta}X$ is a degeneracy of both the nondegenerate simplices $\lambda_{\sigma\theta}(f):\Delta[l]\rightarrow X$ and $\lambda_{\sigma}(\theta(f)):\Delta[l']\rightarrow X$. By the uniqueness in the Eilenberg-Zilberg lemma, these two simplices are equal. Therefore $\lambda_{\sigma\theta}=\lambda_{\sigma}\circ\theta$.

Now let us assume that we have another cocone $\{\mu_\sigma:\pos{\sigma}\rightarrow P,\sigma\in\Delta(X)\}$ where $P$ is a preordered set. We first observe that there is at most one extension $\mu:\pos{X}\rightarrow P$. Given a nondegenerate simplex $(\sigma:\Delta[p]\rightarrow X)\in\pos{X}$, we must have $\mu(\sigma)=\mu(\lambda_{\sigma}([p]=[p]))=\mu_{\sigma}([p]=[p])$. We now prove that this assignment satisfies for all $\sigma\in\Delta(X)$,
\[
\mu\circ\lambda_{\sigma}=\mu_\sigma.
\]

Fix a (possibly degenerate) simplex $\sigma:\Delta[p]\rightarrow X$ and $f:[k]\hookrightarrow [p]\in\pos{\sigma}$. On the one hand,
\[
(\mu\circ\lambda_\sigma)(f)=\mu(ez(\sigma\circ f))=\mu_{ez(\sigma\circ f)}(ez(\sigma\circ f)=ez(\sigma\circ f)).
\]
On the other hand, writing

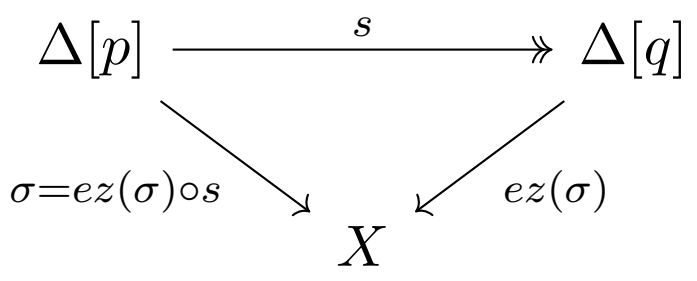
\[
\begin{tikzcd}
\Delta[p]\arrow[rr,"s",two heads]\arrow[dr,"\sigma=ez(\sigma)\circ s"'] & & \Delta[q]\arrow[dl,"ez(\sigma)"] \\
& X &
\end{tikzcd}
\]

for some degeneracy map $s:\Delta[p]\twoheadrightarrow\Delta[q]$, we get

\[
\mu_\sigma(f)=\mu_{ez(\sigma)}(s(f))=\mu_{ez(\sigma)\circ s(f)}\lf((ez(\sigma)\circ s(f))=(ez(\sigma)\circ s(f))\ri).
\]

Now consider the following diagram

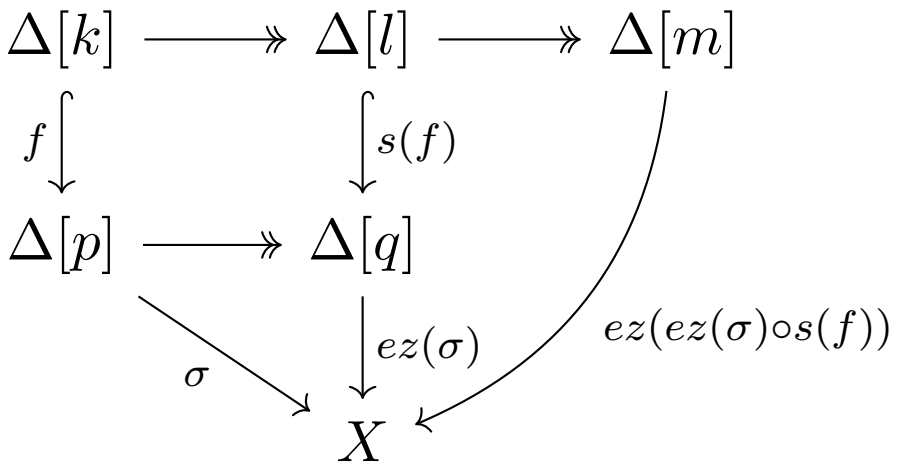
\[
\begin{tikzcd}
\Delta[k]\arrow[r,two heads]\arrow[d,hook,"f"'] & \Delta[l]\arrow[r,two heads]\arrow[d,hook,"s(f)"] & \Delta[m]\arrow[ddl,"ez(ez(\sigma)\circ s(f))",bend left] \\
\Delta[p]\arrow[r,two heads]\arrow[dr,"\sigma"'] & \Delta[q]\arrow[d,"ez(\sigma)"] & \\
& X & 
\end{tikzcd}.
\]

We observe that the simplex $ez(ez(\sigma)\circ s(f)):\Delta[m]\rightarrow X$ is nondegenerate and has for degeneracy the simplex $\sigma\circ f:\Delta[k]\rightarrow X$. Necessarily, we then have
\[
ez(ez(\sigma)\circ s(f))=ez(\sigma\circ f).
\]

Therefore 
\[
(\mu\circ\lambda_\sigma)(f)=\mu_\sigma(f).
\]

It remains to prove that the map $\mu:\pos{X}\rightarrow P$ is compatible with the preorder. Let $\sigma:\Delta[p]\rightarrow X,\tau:\Delta[q]\rightarrow X\in\pos{X}$ such that $\sigma\leq\tau$, that is there exists a face map $f:\Delta[p]\hookrightarrow\Delta[q]$ such that $\tau f=\sigma$.

We compute

\[
\mu(\sigma)=\mu_{\sigma}([p]=[p])=\mu_{\tau}([p]\xrightarrow{f} [q]\leq\mu_{\tau}([q]=[q])=\mu(\tau).
\]

\end{proof}

Endowing a preorder with its Alexandrov topology defines a functor $Alex:\Pre\rightarrow\Top$ which is a left adjoint, see for example \cite[Proposition 2.6]{waas} for an explicit formulation and proof. From the previous lemma and the fact that left adjoints preserve colimits, we deduce the following result. 

\begin{corollary}\label{colim_B(X)}
There is a colimit of topological spaces $\colim{\sigma\in\Delta(X)}{\pos{\sigma}}\cong\pos{X}$ where all posets are endowed with the Alexandrov topology.
\end{corollary}

\begin{rem}
The poset $\pos{X}$ appears as early as \cite{barratt} and has also been used in \cite[Section 2]{waldhausen}. However, in both cases, the focus is on its nerve rather than on the poset itself. In the latter reference, they also use the notation $\pos{X}$ for the nerve of this poset, which they call the "Barratt nerve". To our knowledge, the fact that this poset is obtained as a left Kan extension in the category of topological spaces has not appeared before in the literature.
\end{rem}

\begin{definition}
Let $X$ be a simplicial set. We define the stratification $\pi_X:\norm{X}\rightarrow\pos{X}$ by $\colim{\sigma\in\Delta(X)}{\lf(\Delta^\sigma\rightarrow\pos{\sigma}\ri)}$.
\end{definition}

\begin{rem}
Let $t\in\norm{X}$ and $\sigma\in\Delta(X)$ a non degenerate simplex such that $t$ lies in the interior of the simplex $\Delta^\sigma$ under the quotient map $q_\sigma:\Delta^\sigma\rightarrow\norm{X}$. Consider the diagram

\[
\begin{tikzcd}
\Delta^\sigma\arrow[d,"\pi_\sigma"']\arrow[r,"q_\sigma"] & \norm{X}\arrow[d,"\pi_X"] \\
\pos{\sigma}\arrow[r,"\lambda_\sigma"'] & \pos{\sigma}
\end{tikzcd}
\]
where $\lambda_\sigma:\pos{\sigma}\rightarrow\pos{X}$ is the map defined in the proof of \Cref{colim_preord}. Identifying $t$ with its preimage in $\Delta^\sigma$, we then get
\[
\pi_X(q_\sigma(t))=\lambda_\sigma(\pi_\sigma(t))=\lambda_\sigma(\sigma)=\sigma
\]
since $t$ belongs to the top stratum of $\Delta^\sigma$ and $\sigma$ is non degenerate. This justifies the claim that $\pi_X$ sends a point to the maximal nondegenerate simplex it belongs to. 
\end{rem}

\begin{rem}\label{functoriality_pi}
Since the indexing category $\Delta(X)$ is functorial in $X$, this defines a functor $\pi_{-}:\sSet\rightarrow\text{Arr}(\Top)$.
\end{rem}

\section{Conicality of the stratification}
Let $X$ be a locally finite simplicial set, that is each simplex is a face of a finite number of nondegenerate simplices. In order to define the exit path category $\exit{\pos{X}}{\norm{X}}$, using Lurie's model, we need to prove that the $\pos{X}$-stratification is conical. This section is solely devoted to proving this claim, since it is quite technical.

We start by proving the claim for a finite simplicial set $X$. The locally finite case will follow from this almost directly.

\begin{proposition}\label{conicality}
Let $X$ be a finite simplicial set. Then the stratification $\pi_X:\norm{X}\rightarrow\pos{X}$ is conical.
\end{proposition}

The proof consists of an explicit construction by induction over the filtration of $X$ by its nondegenerate simplices. In the following proof, we will often abuse notation, identifying the image of a topological space $U$ under an embedding with $U$ itself.

Fix an ordering $\sigma_1,\sigma_2,\cdots,\sigma_N$ on the nondegenerate simplices such that the simplices are nondecreasing in dimension. This sequence of simplices induces the filtrations
\begin{align*}
\emptyset:=X^0\hookrightarrow X^1\hookrightarrow\cdots\hookrightarrow X^N=X, \\
\emptyset:=\norm{X^0}\hookrightarrow \norm{X^1}\hookrightarrow\cdots\hookrightarrow \norm{X^N}=\norm{X}
\end{align*}

where $X^j$ is the smallest subsimplicial set containing all simplices $\sigma_1,\cdots,\sigma_j$. Note that both the successive subsimplicial sets $\{X^j\}_{1\leq j\leq N}$ and their geometric realizations are obtained by iterative pushouts. Namely, for each $j$, the map $\Delta[\sigma_j]\rightarrow X$ restricts to a map $\partial\Delta[\sigma_j]\rightarrow X$ which factors through $X^{j-1}$ by construction. Then we have a pushout square of simplicial sets
\[
\begin{tikzcd}
\partial\Delta[\sigma_j]\arrow[r,hook]\arrow[d] & \Delta[\sigma_j]\arrow[d] \\
X^{j-1}\arrow[r] & X^j\arrow[ul,"\ulcorner",phantom,very near start]
\end{tikzcd}
\]

and similarly for $\norm{X^j}$. 

Now fix a point $x\in\norm{X}$ and let $\sigma_k$ be the stratum it lies in. Let $U$ be an open neighbourhood contained in the interior $\overset{\circ}{\Delta}{}^{\sigma_k}$. We may assume that $U$ is locally compact, for example by taking $U$ to be homeomorphic to an open ball. The benefit of this assumption comes from the following result.

\begin{lemma}\cite[Proposition A.14]{hatcher}
Let $U$ be a locally compact topological space. Then the functor $U\times -:\Top\rightarrow\Top$ is a left adjoint.
\end{lemma}

In particular, taking the product with $U$ preserves colimits.

Now for any $j\geq k$, we prove by induction that $x$ admits a conical neighbourhood $\iota^j:U\times C(L^j)\hookrightarrow \norm{X^j}$. To do so, assume that we have a conical neighbourhood $\iota^{j-1}:U\times C(L^{j-1})\hookrightarrow\norm{X^{j-1}}$ around $x$. Taking the quotient map $q^\partial:\partial\Delta^{\sigma_j}\rightarrow\norm{X^{j-1}}$, it restricts to a conical neighbourhood $U\times C(L^\partial)\hookrightarrow q^\partial(\partial\Delta^{\sigma_j})$. The preimage of this conical neighbourhood under $q^\partial$ gives a mapping cylinder $U\times M(f^\partial)\hookrightarrow\partial\Delta^{\sigma_j}$. In order to control this mapping cylinder, one needs mapping cylinders $U\times M(f_\tau)$ on all simplices $\Delta^\tau$ and compatibility between them. In particular, one needs to define such mapping cylinders for degenerate simplices. \Cref{section_degeneracies} addresses this problem. 

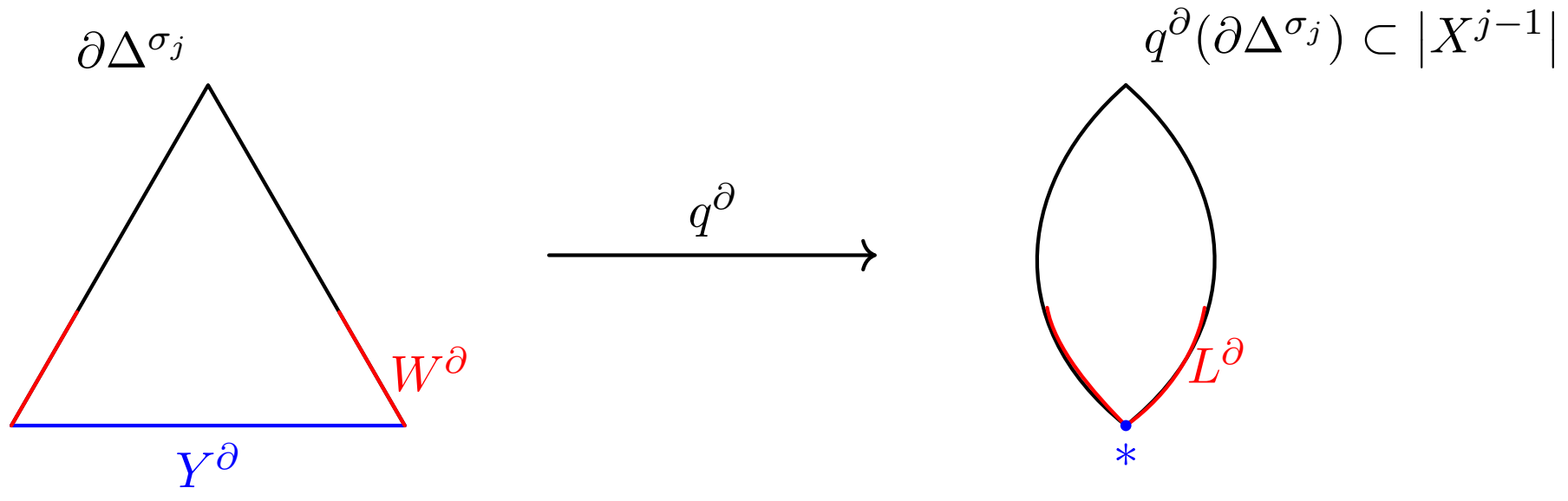
\begin{figure}[ht]
\centering
\begin{tikzpicture}[scale=1, line cap=round, line join=round]

\coordinate (A) at (1,0);
\coordinate (B) at (4,0);
\coordinate (C) at (2.5,{3*sqrt(3)/2});

\coordinate (Ared) at (1.5,{sqrt(3)/2});
\coordinate (Bred) at (3.5,{sqrt(3)/2});

\draw[thick] (A) -- (C) -- (B);
\draw[blue, thick] (A) -- (B);

\draw[red, thick] (A) -- (Ared);
\draw[red, thick] (B) -- (Bred);

\node[above left] at (C) {$\partial\Delta^{\sigma_j}$};
\node[blue, below] at ($(A)!0.5!(B)$) {$Y^\partial$};
\node[red, right] at ($(B)!0.5!(Bred)$) {$W^\partial$};

\draw[->, thick] (5.1,1.3) -- node[above] {$q^\partial$} (7.6,1.3);

\coordinate (P) at (9.5,0);          
\coordinate (T) at (9.5,{3*sqrt(3)/2});

\coordinate (PredL) at (8.90,0.9);
\coordinate (PredR) at (10.10,0.9);

\draw[thick] (P) .. controls (8.5,0.8) and (8.7,1.9) .. (T);
\draw[thick] (P) .. controls (10.5,0.8) and (10.3,1.9) .. (T);

\draw[red, thick] (P) .. controls (9.25,0.25) and (8.95,0.6) .. (PredL);
\draw[red, thick] (P) .. controls (9.85,0.25) and (10.05,0.6) .. (PredR);

\fill[blue] (P) circle (1.2pt);

\node[above right] at (T) {$q^\partial(\partial\Delta^{\sigma_j})\subset\norm{X^{j-1}}$};
\node[blue, below] at (P) {$\ast$};
\node[red, right] at ($(P)!0.55!(PredR)$) {$L^\partial$};

\end{tikzpicture}
\caption{The induced mapping cylinder.}
\label{fig:induced-mapping-cylinder}
\end{figure}

We then extend the conical neighbourhood to $\norm{X^j}$ by taking a collar of $U\times M(f^\partial)$ inside $\Delta^{\sigma_j}$ and attach it to the conical neighbourhood $U\times C(L^{j-1})\hookrightarrow\norm{X^{j-1}}$. This gives the desired conical neighbourhood $U\times C(L^j)\hookrightarrow\norm{X^j}$. In \Cref{section_mapping_cylinders}, we explain how taking a collar of a mapping cylinder gives another mapping cylinder. 

This discussion leads to the stronger inductive hypotheses:

\begin{enumerate}
\item[(IH1)] There exists a stratified space $\pi^j:L^j\rightarrow\pos{X^j}_{\sigma_k<}$ and a conical neighbourhood $\iota^j:U\times C(L^j)\hookrightarrow \norm{X^j}$ where $L^j$ is compact.
\item[(IH2)] There is a functor $f^j:\Delta^{inj}(X^j)\rightarrow\text{Arr}(\Top),(\tau\mapsto f^j_\tau:W_\tau\rightarrow Y_\tau)$, where all topological spaces $W_\tau^j, Y_\tau^j$ are compact.
\item[(IH3)] The functor $f^j$ extends to a functor 
\begin{align*}
& \gamma^j:\Delta^{inj}(X^j)\rightarrow\text{Arr}(\Top)_{/(L^j\rightarrow\ast)}, \\
& \gamma^j_\tau=(f^j_\tau:W_\tau^j\rightarrow Y_\tau^j)\rightarrow(L^j\rightarrow\ast).
\end{align*}
\item[(IH4)] For each simplex $\tau\in\Delta^{inj}(X^j)$, there is a natural transformation $\eta^j:U\times M(f^j)\Rightarrow\Delta^-$ fitting into a pullback square
\[
\begin{tikzcd}
U\times M(f^j_\tau)\arrow[d,"\eta^j_\tau"']\arrow[r,"U\times(\gamma^j_\tau)^\lhd"]\arrow[dr,"\lrcorner",very near start,phantom]  & U\times C(L^j)\arrow[d,"\iota^j"] \\
\Delta^\tau\arrow[r,"q_\tau"'] & \norm{X^j}
\end{tikzcd}.
\]  
\end{enumerate}

From now on, we suppress the stage index on the spaces $W_\tau$ and $Y_\tau$; it will be determined by the functor $f^j$ under consideration.

Before delving into the proof, we establish some notation. We write $0:\ast\rightarrow\R_{\geq 0}$ for the inclusion at $0\in\R_{\geq 0}$. Likewise, $1:\ast\rightarrow [0,1]$ denotes the inclusion at $1$. Moreover, given a map $f:W\rightarrow Y$ and a topological space $Z$, we write $Z\times f:Z\times W\rightarrow Z\times Y$ for the induced map given by the identity on the $Z$-component, and define likewise $f\times Z$.
Given a map $f:W\rightarrow Y$, we write $M(f)$ for the open mapping cylinder defined by the pushout
\[
\begin{tikzcd}
W\arrow[r,"f"]\arrow[d,"W\times 0"'] & Y \arrow[d,"j_f"] \\
W\times \R_{\geq 0}\arrow[r,"q_f"'] & M(f)\arrow[ul,phantom,"\ulcorner",near start] 
\end{tikzcd}
\]
and $\bar{M}(f)$ for the closed mapping cylinder, that is, the same pushout with $W\times \R_{\geq 0}$ replaced by $W\times [0,1]$. We denote its canonical retraction to $Y$ by $r_f:\bar{M}(f)\rightarrow Y$. Note that whenever $f$ is given by a map $L\rightarrow\ast$ with $L$ compact, we recover the cone $C(L)$ by \Cref{teardrop_topology}. All cones occurring throughout the proof are taken from compact topological spaces so this will not cause confusion.

Given maps of topological spaces $Y\xleftarrow{f} W\xrightarrow{g} L$, we define $\bar{D}(f,g)$ to be the double mapping cylinder
\[
\begin{tikzcd}
W\arrow[r,"g"]\arrow[d,"W\times 1"'] & L \arrow[d] \\
\bar{M}(f)\arrow[r] & \bar{D}(f,g)\arrow[ul,phantom,"\ulcorner",near start] 
\end{tikzcd}
\]
where we abuse notation by writing $W\times 1$ for the composition $W\rightarrow W\times [0,1]\rightarrow\bar{M}(f)$.

Also, given a map $g:(f:W\rightarrow Y)\rightarrow (f':W'\rightarrow Y')$, we write $g^\lhd:M(f)\rightarrow M(f')$ for the induced map on mapping cylinders.

Finally, for a simplex $\sigma\in\Delta(X)$, we write $\partial\sigma$ the full subcategory of $\Delta^{inj}(X)_{/\sigma}$ having the same objects except the terminal object $\sigma=\sigma$. In other words, objects are strict face maps $\tau\rightarrow\sigma$ and morphisms are face maps between them.

\subsection{Recognizing mapping cylinders}\label{section_mapping_cylinders}

\begin{lemma}\label{reparametrization_square}
The map $\mu:[0,1]\times \R_{\geq 0}\rightarrow \R_{\geq 0}\times\R_{\geq 0},(r,t)\mapsto (tr,t(1-r))$ exhibits $\R_{\geq 0}\times\R_{\geq 0}$ as the cone $C([0,1])$.
\end{lemma}

\begin{proof}
Since for any $r,r'\in [0,1],\mu(r,0)=\mu(r',0)=(0,0)$, we get an induced map $\psi:C([0,1])\rightarrow \R_{\geq 0}\times\R_{\geq 0}$ sending the tip $\ast$ to $(0,0)$. Conversely, define the map $\varphi:\R_{\geq 0}\times\R_{\geq 0}\rightarrow C([0,1])$ by
\[
\varphi(u,v)=
\lf\{
\begin{array}{ll}
\ast & \text{if }(u,v)=(0,0) \\
\lf(\frac{u}{u+v},u+v\ri) & \text{otherwise}.
\end{array}
\ri.
\]

This map is continuous: it is straightforward outside the basepoint. Moreover given an open neighbourhood $U\subset C([0,1])$ containing $\ast$, it contains an open subset $U_{\epsilon}=\{\ast\}\sqcup\{(r,t)\in U,t<\epsilon\}$ for some $\epsilon>0$. Then $\varphi^{-1}(U)$ contains the subset $\varphi^{-1}(U_\epsilon)=\{(u,v)\in\R_{\geq 0}\times\R_{\geq 0},u+v<\epsilon\}$ which is open and contains $(0,0)$.

We now check that $\psi,\varphi$ are inverses of each other. We have already shown that $\psi(\ast)=(0,0)$ and $\varphi(0,0)=\ast$. For $(r,t)\in C([0,1]),t>0$, observing that $tr+t(1-r)=t$, we have
\[
\varphi(\psi(r,t))=\varphi(tr,t(1-r))=\lf(\frac{tr}{t},t\ri)=(r,t)
\]
and for $(u,v)\in \R_{\geq 0}\times\R_{\geq 0},(u,v)\neq (0,0)$, we have
\[
\psi(\varphi(u,v))=\psi\lf(\frac{u}{u+v},u+v\ri)=\lf(u,\frac{v}{u+v}(u+v)\ri)=(u,v).
\]

\end{proof}

For a continuous map $f:W\rightarrow Y$, consider the topological space $M(f)\times\R_{\geq 0}$ pictured in the figure below. Each point on the vertical axis is identified with a copy of $Y$ and each point in the interior or on the horizontal axis is identified with a copy of $W$. The $\R_{\geq 0}$-coordinate is determined by the vertical axis.

\[
\begin{tikzpicture}[scale=1.2, line cap=round, line join=round]
  \draw[thick] (0,0) -- (4,0);
  \draw[thick] (0,0) -- (0,4);

  \node[right] at (4,0) {$\cdots$};
  \node[above] at (0,4) {$\cdots$};

  \node[below left] at (0,0) {$Y$};
  \node[left] at (0,2) {$Y$};
  \node[red,below] at (2.5,0) {$W\times \{1\}$};
  \node at (2.2,2.4) {$W$};

  \draw[red, thick] (0,2.6) -- (2.6,0);
  \node[red, above right] at (1.15,1.15) {$\bar{M}(f)$};
\end{tikzpicture}
\]

Following \Cref{reparametrization_square}, the diagonal direction suggests it can be identified as the mapping cylinder of the retraction $r_f:\bar{M}(f)\rightarrow Y$. A slice of this mapping cylinder is depicted in red. With this identification, the inclusion $W\xrightarrow{W\times 1}\bar{M}(f)$ induces a map $M(f)\rightarrow M(r_f)$ whose image is identified with the lower horizontal axis. The following lemma makes these ideas precise.

\begin{lemma}\label{cone_prism}
Let $f:W\rightarrow Y$ be a continuous map. Then there is a homeomorphism $h:M(r_f)\xrightarrow{\cong} M(f)\times \R_{\geq 0}$ fitting into a commutative triangle
\[
\begin{tikzcd}
& M(f)\arrow[dl,"(W\times 1)^\lhd"']\arrow[dr,"M(f)\times 0"] & \\
M(r_f)\arrow[rr,"h"'] & & M(f)\times\R_{\geq 0}
\end{tikzcd}
\]
\end{lemma}

\begin{proof} 

Consider the diagram
\[
\resizebox{\textwidth}{!}{%
\begin{tikzcd}[ampersand replacement=\&]
W\times 0\times \R_{\geq 0}
  \arrow[dr,hook]
  \arrow[dd,"W\times \mu"',near start]
  \arrow[rr,"f\times\R_{\geq 0}",near start]
\& \&
Y\times 0\times \R_{\geq 0}
  \arrow[dr,hook,"\bar{j}_f\times\R_{\geq 0}"]
  \arrow[dd,"Y\times\mu",near start]
\&
\\
\&
W\times [0,1]\times\R_{\geq 0}
  \arrow[rr,"\bar{q}_f",near start,crossing over]
  \arrow[dd,"W\times\mu"',near start,crossing over]
\& \&
\bar{M}(f)\times \R_{\geq 0}
  \arrow[dd,"\mu_f",dotted,near start]
\\
W\times 0 \times\R_{\geq 0}
  \arrow[rr]
  \arrow[dr,hook]
\& \&
Y\times 0\times\R_{\geq 0}
  \arrow[dr,hook,"j_f\times\R_{\geq 0}"]
\&
\\
\&
W\times \R_{\geq 0}\times\R_{\geq 0}
  \arrow[rr,"q_f"',near start]
\& \&
M(f)\times\R_{\geq 0}
\end{tikzcd}%
}
\]

where the top and bottom faces are pushouts. By definition of $\mu$, this diagram commutes, hence induces a map $\mu_f:\bar{M}(f)\times\R_{\geq 0}\rightarrow M(f)\times\R_{\geq 0}$ fitting into a commutative square

\[
\begin{tikzcd}
\bar{M}(f)\times 0\arrow[r,"r_f"]\arrow[d,hook] & Y\times 0\arrow[d,hook,"j_f\times 0"] \\
\bar{M}(f)\times\R_{\geq 0}\arrow[r,"\mu_f"'] & M(f)\times\R_{\geq 0}
\end{tikzcd}.
\]
Moreover $\mu_f$ is surjective since the maps 
\begin{align*}
(\bar{q}\circ(W\times\mu):W\times [0,1]\times\R_{\geq 0}\rightarrow M(f)\times\R_{\geq 0}, \\
(\bar{\iota}\times\R_{\geq 0})\circ (Y\times\mu):Y\times 0\times\R_{\geq 0}\rightarrow M(f)\times\R_{\geq 0})
\end{align*}
are jointly surjective.

Now consider a commutative square
\[
\begin{tikzcd}
\bar{M}(f)\times 0\arrow[r,"r_f"]\arrow[d,hook,"\bar{M}(f)\times 0"'] & Y\arrow[d,"\alpha"] \\
\bar{M}(f)\times\R_{\geq 0}\arrow[r,"\beta"'] & T
\end{tikzcd}.
\]

Observe that given an element $x\in\bar{M}(f)$, we have $\beta(x,0)=\alpha(r_f(x))$ which is independent of the radial direction of $\bar{M}(f)$. Because of this, \Cref{reparametrization_square} provides a unique map $\bar{\beta}:W\times\R_{\geq 0}\times\R_{\geq 0}\rightarrow T$ satisfying $\bar{\beta}\circ (W\times\mu)=\beta\circ q$. We then get a commutative square

\[
\begin{tikzcd}
W\times 0\times\R_{\geq 0}\arrow[r]\arrow[d] & Y\times 0\times\R_{\geq 0}\arrow[d,"\beta\circ(\iota\times\R_{\geq 0})"] \\
W\times\R_{\geq 0}\times\R_{\geq 0}\arrow[r,"\bar{\beta}"'] & T
\end{tikzcd}
\]
inducing a unique map $\xi:M(f)\times\R_{\geq 0}\rightarrow T$ such that $\xi\circ(\bar{\iota}\times\R_{\geq 0})=\beta\circ(\iota\times\R_{\geq 0}),\xi\circ\bar{q}=\bar{\beta}$.

On the one hand, we observe that both $\xi\circ\mu_f,\beta$ satisfy the universal property applied to the upper pushout
\[
\begin{tikzcd}
W\times 0\times\R_{\geq 0}\arrow[r,"f\times\R_{\geq 0}"]\arrow[d,hook] & Y\times\{0\}\times\R_{\geq 0}\arrow[d,hook] \\
W\times [0,1]\times\R_{\geq 0}\arrow[r,"q"'] & \bar{M}(f)\times\R_{\geq 0}
\end{tikzcd}
\]

hence $\xi\circ\mu_f=\beta$. In particular, we get the following commutative diagram.

\[
\begin{tikzcd}
Y\times 0\times 0\arrow[d,"\iota\times 0"']\arrow[r,"\bar{\iota}\times 0"] & M(f)\times\R_{\geq 0}\arrow[d,"\xi"] \\
\bar{M}(f)\times\R_{\geq 0}\arrow[r,"\beta"']\arrow[ur,"\mu_f"'] & T
\end{tikzcd}
\]

On the other hand, the identity $\beta\circ(\iota\times 0)=\alpha$ provides $\xi\circ(\bar{\iota}\times 0)=\beta\circ(\iota\times 0)=\alpha$.
We have thus proved the existence of a map $\xi:M(f)\times\R_{\geq 0}\rightarrow T$ satisfying $\xi\circ(\bar{\iota}\times 0)=\alpha,\xi\circ\mu_f=\beta$. Moreover, if $\xi':M(f)\times\R_{\geq 0}\rightarrow T$ is another map satisfying these identities, then the surjectivity of $\mu_f$ implies $\xi=\xi'$. Thus such an extension $\xi$ is unique, thereby showing that the square 

\[
\begin{tikzcd}
\bar{M}(f)\times 0\arrow[r,"r_f"]\arrow[d,hook,"\bar{M}(f)\times 0"'] & Y\times 0\arrow[d,hook,"j_f\times 0"] \\
\bar{M}(f)\times\R_{\geq 0}\arrow[r,"\mu_f"'] & M(f)\times\R_{\geq 0}
\end{tikzcd}
\]
is a pushout.

We then get a commutative diagram
\[
\begin{tikzcd}
\bar{M}(f)\times 0\arrow[d,hook,"\bar{M}(f)\times 0"']\arrow[r,"r_f"] & Y\times 0\arrow[ddr,"j_f\times 0",bend left]\arrow[d,"j_{r_f}",hook] & \\
\bar{M}(f)\times\R_{\geq 0}\arrow[drr,bend right,"\mu_f"']\arrow[r,"q_{r_f}"'] & M(r_f)\arrow[dr,"h",dotted] & \\
& & M(f)\times\R_{\geq 0}
\end{tikzcd}
\]
where the upper-left square is a pushout, hence gives the desired homeomorphism 
\[
h:M(r_f)\xrightarrow{\cong}M(f)\times\R_{\geq 0}.
\]
Now consider the diagram
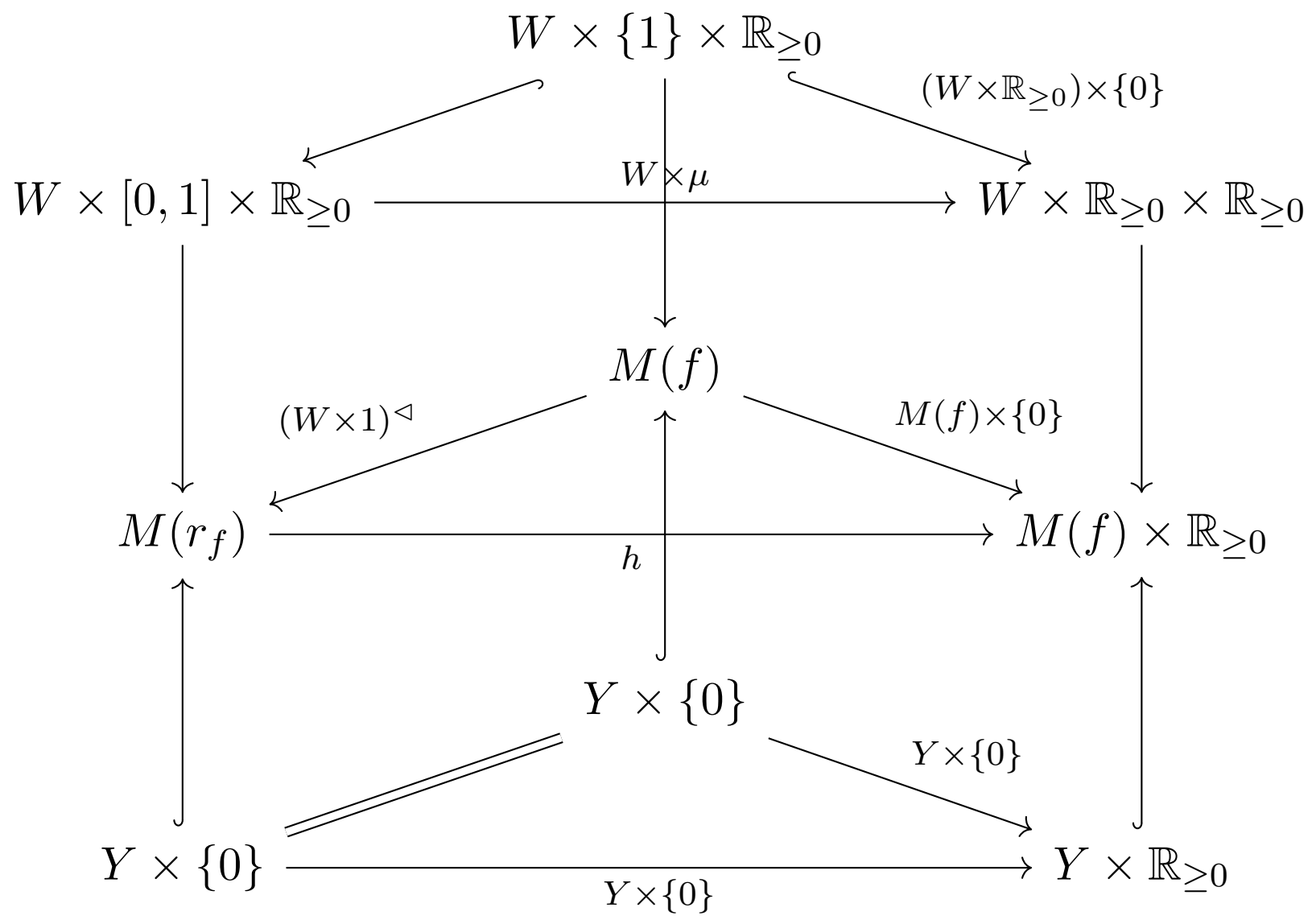
\[
\begin{tikzcd}
& W\times\{1\}\times\R_{\geq 0}\arrow[dl,hook]\arrow[dd]\arrow[dr,hook,"(W\times\R_{\geq 0})\times\{0\}"] & \\
W\times [0,1]\times\R_{\geq 0}\arrow[dd]\arrow[rr,"W\times\mu"] & & W\times\R_{\geq 0}\times\R_{\geq 0}\arrow[dd] \\
& M(f)\arrow[dl,"(W\times 1)^\lhd"']\arrow[dr,"M(f)\times\{0\}"] & \\
M(r_f)\arrow[rr,"h"'] & & M(f)\times\R_{\geq 0} \\
& Y\times\{0\}\arrow[dl,equal]\arrow[dr,"Y\times\{0\}"]\arrow[uu,hook] & \\
Y\times\{0\}\arrow[rr,"Y\times\{0\}"']\arrow[uu,hook] & & Y\times\R_{\geq 0}\arrow[uu,hook] 
\end{tikzcd}
\]
where all lateral rectangles commute by construction and so does the lower triangle.
The upper triangle also commutes because for $(w,1,t)\in W\times\{1\}\times\R_{\geq 0}$, we have $(W\times\mu)(w,1,t)=(w,t,0)$ hence the top right map is indeed identified with the inclusion $(W\times\R_{\geq 0})\times\{0\}:W\times\R_{\geq 0}\hookrightarrow W\times\R_{\geq 0}\times\R_{\geq 0}$. By joint surjectivity of the vertical maps, the middle triangle also commutes, as desired.

\end{proof}

\begin{lemma}\label{prism_lemma}
Consider a commutative square of maps between compact topological spaces
\[
\begin{tikzcd}
W\arrow[r,"g"]\arrow[d,"f"'] & L\arrow[d] \\
Y\arrow[r] & \ast
\end{tikzcd}
\]
where $L$ is compact. Then there is a pushout square

\[
\begin{tikzcd}
M(f)\arrow[r,"g^\lhd"]\arrow[d,"(W\times 1)^\lhd"'] & C(L)\arrow[d] \\
M(r_f)\arrow[r] & C(\bar{D}(f,g))\arrow[ul,"\ulcorner",phantom,very near start]
\end{tikzcd}.
\]

\end{lemma}

We first give a visual explanation of this result. As for \Cref{cone_prism}, the pushout can be pictured as follows.

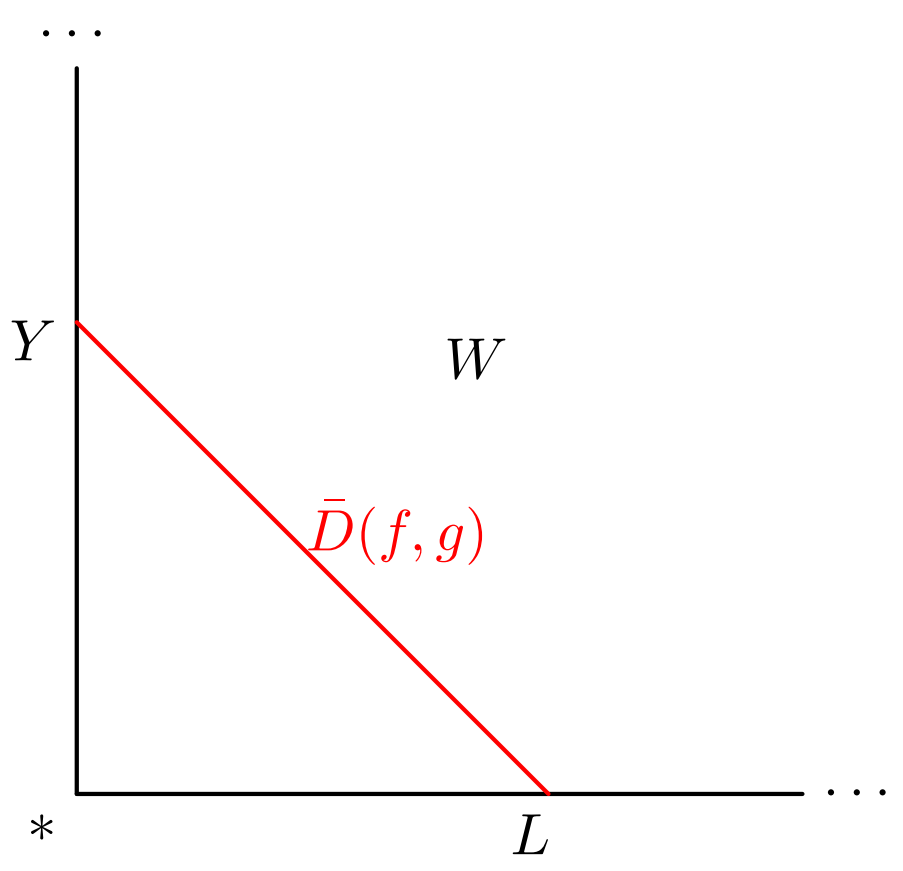
\[
\begin{tikzpicture}[scale=1.2, line cap=round, line join=round]
  \draw[thick] (0,0) -- (4,0);
  \draw[thick] (0,0) -- (0,4);

  \node[right] at (4,0) {$\cdots$};
  \node[above] at (0,4) {$\cdots$};

  \node[below left] at (0,0) {$\ast$};
  \node[below] at (2.5,0) {$L$};
  \node[left] at (0,2.5) {$Y$};
  \node at (2.2,2.4) {$W$};

  \draw[red, thick] (0,2.6) -- (2.6,0);
  \node[red, above right] at (1.15,1.15) {$\bar{D}(f,g)$};
\end{tikzpicture}
\]

Again, the diagonal direction suggests this picture can be identified with the cone $C(\bar{D}(f,g))$.

\begin{proof}
First observe that $\bar{D}(f,g)$ is compact since $W,L,Y$ are compact. Thus there is no ambiguity in writing $C(\bar{D}(f,g))$ as the usual pushout.

We first exhibit the maps $g^\lhd,(W\times 1)^\lhd$. The mapping cylinder construction is functorial in maps between topological spaces. Hence, the commutative square $g:(f:W\rightarrow Y)\rightarrow (L\rightarrow\ast)$ induces a map $g^\lhd:M(f)\rightarrow C(L)$. Likewise, the inclusion $W\times\{1\}\hookrightarrow W\times [0,1]$ induces an inclusion $W\hookrightarrow\bar{M}(f)$ which fits into a commutative triangle
\[
\begin{tikzcd}
W\arrow[dr,"f"']\arrow[rr,"W\times 1"] & & \bar{M}(f)\arrow[dl,"r_f"] \\
& Y & 
\end{tikzcd}.
\] 
This provides the map $(W\times 1)^\lhd:M(f)\rightarrow M(r_f)$.

Now consider a commutative square 

\[
\begin{tikzcd}
M(f)\arrow[r,"g^\lhd"]\arrow[d,"(W\times 1)^\lhd"'] & C(L)\arrow[d,"\alpha"] \\
M(r_f)\arrow[r,"\beta"'] & T
\end{tikzcd},
\]
and the diagram

\[
\begin{tikzcd}
W\times \R_{\geq 0}\arrow[d,"W\times 1\times \R_{\geq 0}"']\arrow[r,"g\times \R_{\geq 0}"] & L\times \R_{\geq 0}\arrow[r,"q_L"]\arrow[d] & C(L)\arrow[dd,"\alpha"] \\
\bar{M}(f)\times \R_{\geq 0}\arrow[d,"q_{r_f}"']\arrow[r] & \bar{D}(f,g)\times \R_{\geq 0} & \\
M(r_f)\arrow[rr,"\beta"'] & & T
\end{tikzcd}
\]

where the map $W\times 1\times\R_{\geq 0}:W\times\R_{\geq 0}\rightarrow\bar{M}(f)\times\R_{\geq 0}$ is induced by the composition $W\times\{1\}\rightarrow W\times [0,1]\rightarrow\bar{M}(f)$. We observe that the top horizontal composition
\[
W\times\R_{\geq 0}\xrightarrow{g\times\R_{\geq 0}}L\times \R_{\geq 0}\xrightarrow{q_L} C(L)
\] 
is equal to the composition $W\times\R_{\geq 0}\rightarrow M(f)\xrightarrow{g^\lhd}C(L)$. Likewise, the left vertical composition
\[
W\times\R_{\geq 0}\xrightarrow{W\times 1\times\R_{\geq 0}}\bar{M}(f)\times\R_{\geq 0}\xrightarrow{q_{r_f}} M(r_f)
\]
is equal to the composition 
\[
W\times\R_{\geq 0}\rightarrow M(f)\xrightarrow{(W\times 1)^\lhd}M(r_f).
\]

Combining with the commutativity of the diagram

\[
\begin{tikzcd}
W\times\R_{\geq 0}\arrow[r] & M(f)\arrow[r,"g^\lhd"]\arrow[d,"(W\times 1)^\lhd"'] & C(L)\arrow[d,"\alpha"] \\
& M(r_f)\arrow[r,"\beta"'] & T
\end{tikzcd}
\]

implies the commutativity of the whole square. Moreover, the upper-left square is a pushout by definition and hence commutes. Thus the whole diagram commutes.

In addition, since the upper-left square is a pushout, it induces a unique map $\psi:\bar{D}(f,g)\times \R_{\geq 0}\rightarrow T$. For any $x\in\bar{D}(f,g)$, we also have $\psi(x,0)=\alpha(\ast)$, thereby providing a unique map $\bar{\psi}:C(\bar{D}(f,g))\rightarrow T$ fitting into a commutative diagram

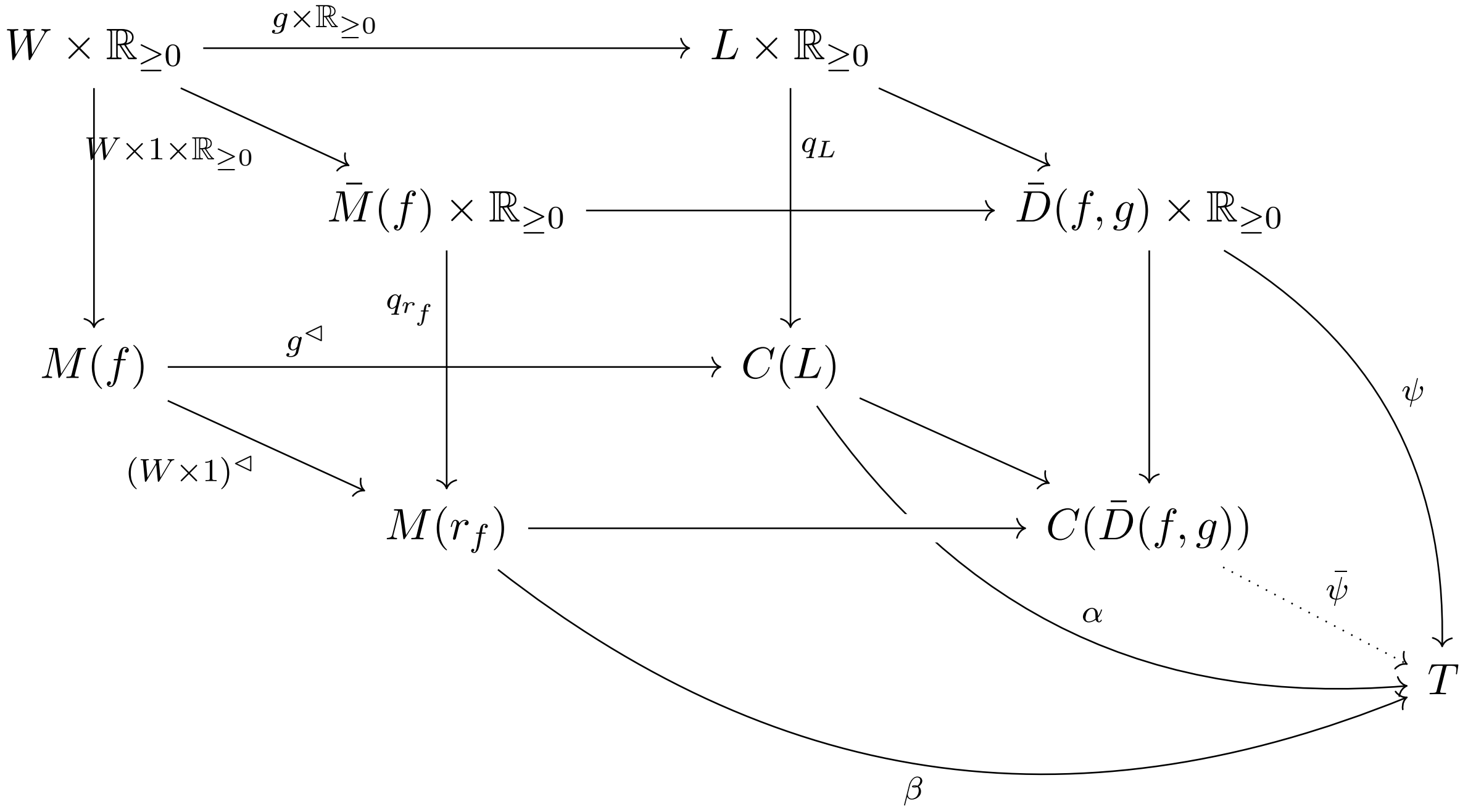
\[
\resizebox{\textwidth}{!}{%
\begin{tikzcd}[ampersand replacement=\&]
W\times\R_{\geq 0}
  \arrow[dr,"W\times 1\times\R_{\geq 0}"']
  \arrow[rr,"g\times\R_{\geq 0}",near start]
  \arrow[dd]
\& \&
L\times\R_{\geq 0}
  \arrow[dr]
  \arrow[dd,"q_L",near start]
\& \&
\\
\&
\bar{M}(f)\times\R_{\geq 0}
  \arrow[rr]
  \arrow[dd,"q_{r_f}"',near start,crossing over]
\& \&
\bar{D}(f,g)\times\R_{\geq 0}
  \arrow[dd]
  \arrow[dddr,bend left,"\psi"]
\&
\\
M(f)
  \arrow[rr,"g^\lhd",near start]
  \arrow[dr,"(W\times 1)^\lhd"']
\& \&
C(L)
  \arrow[dr]
  \arrow[ddrr,bend right,"\alpha"]
\& \&
\\
\&
M(r_f)
  \arrow[rr,crossing over]
  \arrow[rrrd,bend right,"\beta"']
\& \&
C(\bar{D}(f,g))
  \arrow[dr,dotted,"\bar{\psi}"]
\&
\\
\& \& \& \& T
\end{tikzcd}%
}.
\]

Moreover, if another map $\xi:C(\bar{D}(f,g))\rightarrow T$ fits into the lower square, then necessarily the composition $\bar{D}(f,g)\times\R_{\geq 0}\rightarrow C(\bar{D}(f,g))\xrightarrow{\xi} T$ must be equal to $\psi$. Since the first map is surjective, we necessarily get $\bar{\psi}=\xi$. To summarize, we have proved the existence and uniqueness of a map $\bar{\psi}:C(\bar{D}(f,g))\rightarrow T$ fitting into a commutative diagram
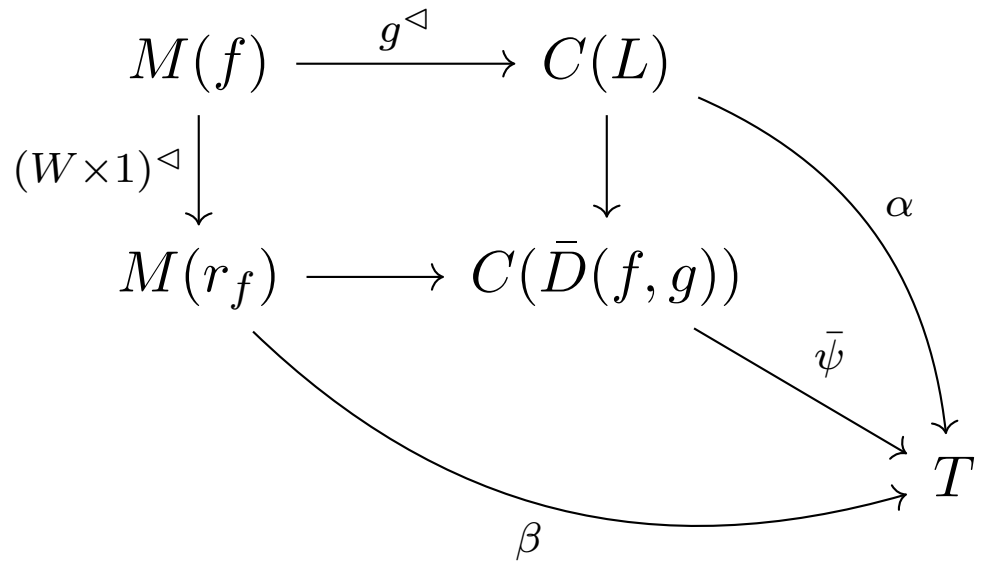
\[
\begin{tikzcd}
M(f)\arrow[d,"(W\times 1)^\lhd"']\arrow[r,"g^\lhd"] & C(L)\arrow[d]\arrow[ddr,bend left,"\alpha"] & \\
M(r_f)\arrow[r]\arrow[drr,bend right,"\beta"'] & C(\bar{D}(f,g))\arrow[dr,"\bar{\psi}"] & \\
& & T
\end{tikzcd}
\]
Hence, the upper-left square is a pushout. 
\end{proof}

\subsection{Lifting conical neighbourhoods through degeneracies}\label{section_degeneracies}

Let $s:\Delta^q\rightarrow\Delta^p$ be a degeneracy map. For two elements $x_0,x_1\in\Delta^p$ with $x_1$ belonging to a lower stratum than $x_0$, we want to define a continuous map $s^{-1}(x_0)\rightarrow s^{-1}(x_1)$ which is also continuous in $x_0,x_1$. The following lemma defines such a map.  

\begin{lemma}\label{translating_fibers}
Let $s:\Delta^q\rightarrow\Delta^p$ be a degeneracy map and $D\subset\Delta^q\times\Delta^p$ be the subspace of pairs $(t,x)$ where $\pi_p(x)\leq\pi_p(s(t))$. Then the assignment $(t,x)\mapsto\rho(t,x)$ with
\[
\rho(t,x)_j=
\lf\{
\begin{array}{ll}
\frac{x_{s(j)} t_j}{s(t)_{s(j)}} & s(t)_{s(j)}\neq 0 \\
0 & \text{otherwise}
\end{array}
\ri.
\]
defines a continuous map $\rho:D\rightarrow\Delta^q$ fitting into a commutative triangle

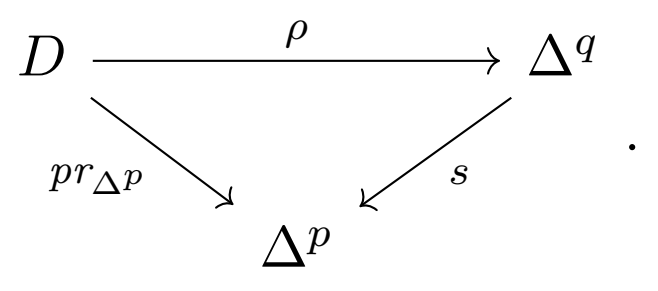
\[
\begin{tikzcd}
D\arrow[rr,"\rho"]\arrow[dr,"\text{pr}_{\Delta^p}"'] & & \Delta^q\arrow[dl,"s"] \\
& \Delta^p & 
\end{tikzcd}.
\]

Moreover, for any $t\in\Delta^q$ and $x,y\in\Delta^p$ such that $\pi_p(y)\leq\pi_p(x)\leq\pi_p(s(t))$, we have
\begin{enumerate}
\item $\rho(t,s(t))=t$,
\item $\rho(\rho(t,x),y)=\rho(t,y)$,
\item if $\pi_p(y)=\pi_p(x)$, $\rho$ restricts to a homeomorphism $\rho(-,y):s^{-1}(x)\rightarrow s^{-1}(y)$ with inverse $\rho(-,x):s^{-1}(y)\rightarrow s^{-1}(x)$.
\end{enumerate}

\end{lemma}

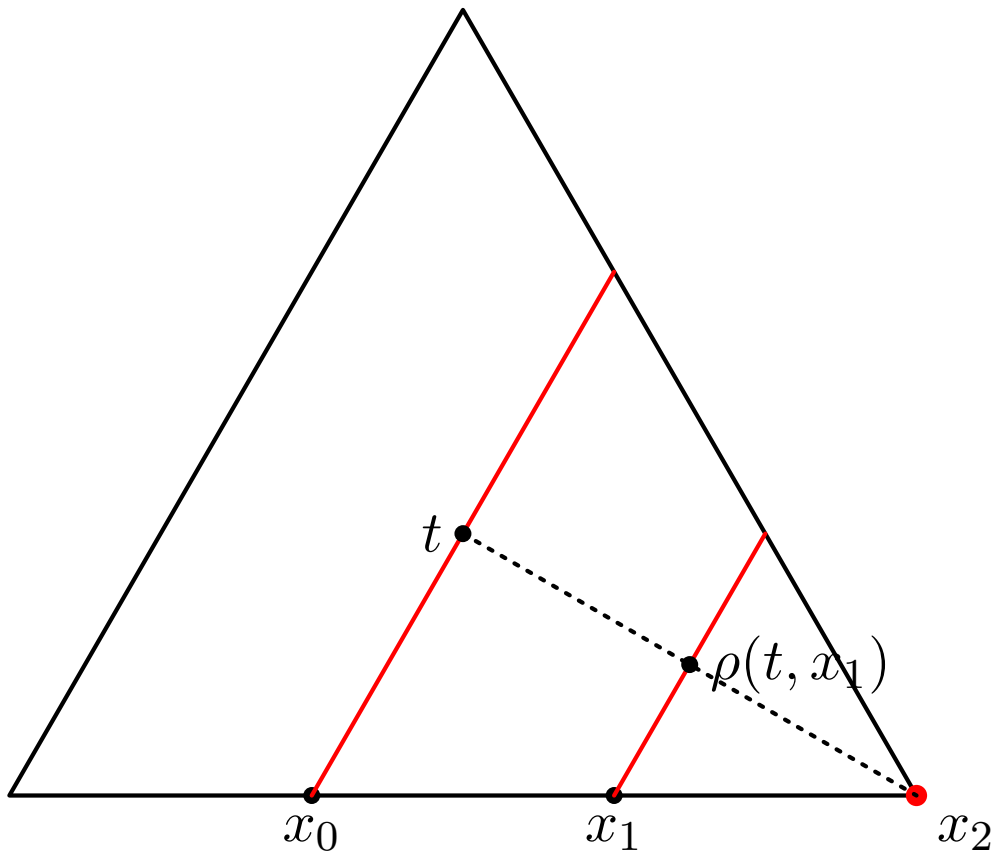
\begin{figure}[ht]
\centering
\begin{tikzpicture}[scale=2, line cap=round, line join=round]
  \coordinate (A) at (0,0);
  \coordinate (B) at (3,0);
  \coordinate (C) at (1.5,{1.5*sqrt(3)});
  \draw[thick] (A) -- (B) -- (C) -- cycle;

  \coordinate (x0) at (1,0);
  \coordinate (x1) at (2,0);
  \coordinate (x2) at (3,0);

  \fill (x0) circle (0.8pt);
  \fill (x1) circle (0.8pt);
  \fill (x2) circle (0.8pt);

  \node[below] at (x0) {$x_0$};
  \node[below] at (x1) {$x_1$};
  \node[below right] at (x2) {$x_2$};

  \coordinate (y0) at (2,{sqrt(3)});
  \coordinate (y1) at (2.5,{0.5*sqrt(3)});

  \draw[red, thick] (x0) -- (y0);
  \draw[red, thick] (x1) -- (y1);
  \fill[red] (x2) circle (1pt);

  \coordinate (t) at (1.5,{0.5*sqrt(3)});
  \fill (t) circle (0.8pt);
  \node[left] at (t) {$t$};

  \draw[dotted, thick] (t) -- (x2);

  \coordinate (rho) at (2.25,{0.25*sqrt(3)});
  \fill (rho) circle (0.8pt);
  \node[right] at (rho) {$\rho(t,x_1)$};
\end{tikzpicture}
\caption{Depicting the map $\rho$.}
\label{fig:rho-map}
\end{figure}

\begin{proof}
By definition, for $(t,x)\in D$ and $j\in [q]$, the component $\rho(t,x)_j$ is well-defined. Moreover
\begin{align*}
\sum_{j\in [q]}\rho(t,x)_j & = \sum_{j\in [q],\rho(t,x)_j\neq 0}\rho(t,x)_j \\
&= \sum_{j\in [q],s(t)_{s(j)}\neq 0}\frac{x_{s(j)}t_j}{s(t)_{s(j)}} \\
&= \sum_{i\in [p],s(t)_i\neq 0}\frac{x_i}{s(t)_i}\sum_{j\in s^{-1}(i)}t_j \\
& = \sum_{i\in [p],s(t)_i\neq 0}x_i \\
& = 1.
\end{align*}
The last equality comes from the fact that $\pi_p(x)\leq\pi_p(s(t))$ so that $x_i\neq 0$ implies $s(t)_i\neq 0$. Also, each component is nonnegative so that $\rho(t,x)\in\Delta^q$. Thus $\rho$ is well-defined.  

We now check the continuity. First observe, from the inequality $t_j\leq s(t)_{s(j)}$, that for any $(t,x)\in D,\rho(t,x)_j\leq x_{s(j)}$. Now let $(t,x)\in D$. We distinguish two cases:
\begin{enumerate}
\item $s(t)_{s(j)}\neq 0$: in this case, there exists an open neighbourhood $U\subset D$ around $(t,x)$ such that for any $(t',x')\in U,s(t')_{s(j)}\neq 0$. Then $(\rho_{|U})_j$ is continuous by continuity of the function $(t',x')\mapsto\frac{x'_{s(j)} t'_j}{s(t')_{s(j)}}$.
\item $s(t)_{s(j)} =0$: in this case, because $\pi_p(x)\leq\pi_p(s(t))$, we also have $x_{s(j)}=0$. Hence we can choose an open neighbourhood $U$ around $(t,x)$ so that for any $(t',x')\in U,x'_{s(j)}$ can be made arbitrarily small. By the previous inequality, $\rho(t',x')_j$ can therefore be made arbitrarily small, which proves continuity.
\end{enumerate}

Now, for $(t,x)\in D$ and $i\in [p]$, we get $s(\rho(t,x))_i=0=x_i$ if $s(t)_{i}=0$ and otherwise

\begin{align*}
s(\rho(t,x))_i &= \sum_{j\in s^{-1}(i)}\rho(t,x)_j \\
&= \frac{x_i}{s(t)_i}\sum_{j\in s^{-1}(i)}t_j \\
&= \frac{x_i}{s(t)_i}s(t)_i \\
&= x_i.
\end{align*}
Therefore $s\circ\rho = \text{pr}_{\Delta^p}$.

Furthermore, the identity $\rho(t,s(t))=t$ is straightforward from the definition. For the other identity, for $j\in [q]$, both $\rho(\rho(t,x),y)_j,\rho(t,y)_j$ are zero exactly when $y_{s(j)}=0$ or $t_j=0$. Otherwise,
\[
\rho(\rho(t,x),y)_j=\frac{y_{s(j)}\rho(t,x)_j}{s(\rho(t,x))_{s(j)}}=\frac{y_{s(j)}}{x_{s(j)}}\frac{x_{s(j)}t_j}{s(t)_{s(j)}}=\rho(t,y)_j.
\]

Finally, for $x,y$ in the same stratum, the identity $s\circ\rho=\text{pr}_{\Delta^p}$ ensures that the restriction $\rho(-,y):s^{-1}(x)\rightarrow\Delta^q$ has image contained in $s^{-1}(y)$, thereby providing a map $\rho(-,y):s^{-1}(x)\rightarrow s^{-1}(y)$. Similarly, we get a map $\rho(-,x):s^{-1}(y)\rightarrow s^{-1}(x)$. Using the two previous identities, for $t\in s^{-1}(y)$, we get
\[
\rho(\rho(t,x),y)=\rho(t,y)=\rho(t,s(t))=t
\]
and similarly for $t\in s^{-1}(x),\rho(\rho(t,y),x)=t$. This proves that these two maps are indeed inverse of each other, hence each one is a homeomorphism.

\end{proof}

Let $f:W\rightarrow Y$ be a continuous map. For $p\geq 0$ and an element $a\in\pos{p}$, consider two stratifications $W\rightarrow\pos{p}_{a<},Y\rightarrow\{a\}\hookrightarrow\pos{p}_{a\leq}$. As for the conical case, this provides a stratification $U\times M(f)\rightarrow\pos{p}_{a\leq}$ whose $a$-stratum is precisely the base $U\times Y$.  

\begin{definition}\label{stratified_cylinder}
A \emph{stratified open embedding} $\eta:U\times M(f)\hookrightarrow\Delta^p$ is an open embedding fitting into a commutative diagram
\[
\begin{tikzcd}
U\times M(f)\arrow[d]\arrow[r,hook,"\eta"] & \Delta^p\arrow[d,"\pi_p"] \\
\pos{p}_{a\leq}\arrow[r,hook] & \pos{p} 
\end{tikzcd}.
\]
\end{definition}

\begin{lemma}\label{lifting_degeneracies}
Let $f:W\rightarrow Y$ be a continuous map between compact stratified spaces as in \Cref{stratified_cylinder} and $\eta:U\times M(f)\hookrightarrow\Delta^p$ a stratified open embedding. For any degeneracy $s:\Delta^q\rightarrow\Delta^p$, there exists a map $f^s:W^s\rightarrow Y^s$ where both spaces are compact and an open embedding $\eta^s:U\times M(f^s)\hookrightarrow\Delta^q$ fitting into a pullback square
\[
\begin{tikzcd}
U\times M(f^s)\arrow[d,"U\times s^\lhd"']\arrow[r,"\eta^s"]\arrow[dr,"\lrcorner",very near start,phantom] & \Delta^q\arrow[d,"s"] \\
U\times M(f)\arrow[r,"\eta"'] & \Delta^p
\end{tikzcd}.
\]

\end{lemma}

\begin{proof}
Fix a basepoint $u\in U$ and define 

\begin{align*}
W^s & :=s^{-1}(\eta(\{u\}\times W\times\{1\})), \\
Y^s & :=s^{-1}(\eta(\{u\}\times Y)).
\end{align*}

Note that both $W^s,Y^s$ are compact since $W,Y$ are compact and $s$ is a proper map.

The map $f^s:W^s\rightarrow Y^s$ is defined by $f^s(t)=\rho(t,f(s(t)))$. This is well-defined since $s(f^s(t))=s(\rho(t,f(s(t))))=f(s(t))\in \eta(Y)$.

Similarly, we define the maps

\begin{align*}
j:U\times Y^s\rightarrow\Delta^q &, j(v,t)=\rho(t,\eta(v,s(t))),\\
F:U\times W^s\times \R_{\geq 0}\rightarrow\Delta^q &, (v,t,r)\mapsto\rho(t,\eta(v,s(t),r))
\end{align*}
where, by abuse of notation, we identify $(v,s(t),r)\in U\times W\times \R_{\geq 0}$ with its quotient image in $U\times M(f)$. With this identification, we get

\begin{align*}
j(v,f^s(t)) &= j(v,\rho(t,f(s(t)))) \\
&=\rho\lf(\rho(t,f(s(t))),\eta(v,f(s(t)))\ri) \\
&= \rho\lf(t,\eta(v,f(s(t)))\ri) \\
&= \rho(t,\eta(v,s(t),0)) \\
&= F(v,t,0).
\end{align*}

Thus $F$ fits into a commutative square

\[
\begin{tikzcd}
U\times W^s\arrow[d,"U\times W^s\times 0"']\arrow[r,"U\times f^s"] & U\times Y^s\arrow[d,"j"] \\
U\times W^s\times \R_{\geq 0}\arrow[r,"F"'] & \Delta^q
\end{tikzcd}.
\]

From this, we obtain a map $\eta^s:U\times M(f^s)\rightarrow\Delta^q$. We show that $\eta^s$ is injective. Assume that $(v,t,r),(v',t',r')\in U\times W^s\times \R_{\geq 0}$ satisfy $F(v,t,r)=F(v',t',r')$. First, looking at the strata they live in, we observe that $r=0$ if and only if $r'=0$. Hence we get two cases:
\begin{enumerate}
\item $r,r'\neq 0$: we have by definition 
\[
\rho(t,\eta(v,s(t),r))=\rho(t',\eta(v',s(t'),r')).
\]
Composing with $s$, we get $\eta(v,s(t),r)=\eta(v',s(t'),r')$. Since $\eta$ is an open embedding, we get $(v,s(t),r)=(v',s(t'),r')$. Now, because $x:=\eta(v,s(t),r)$ and $s(t)$ are in the same stratum, we get by \Cref{translating_fibers} a homeomorphism $\rho(-,x):s^{-1}(s(t))\rightarrow s^{-1}(x)$. In particular, we must have $t=t'$. 
\item $r=r'=0$: it is then a consequence of the injectivity of the map $j:U\times Y^s\rightarrow\Delta^q$, which is obtained by the same argument. 
\end{enumerate}

Moreover, because of the identity $s\circ\rho (t,f(s(t)))=f(s(t))$, we get a commutative square

\[
\begin{tikzcd}
W^s\arrow[d,"s"']\arrow[r,"f^s"] & Y^s\arrow[d,"s"] \\
W\arrow[r,"f"'] & Y
\end{tikzcd}
\]

inducing a map $U\times s^\lhd:U\times M(f^s)\rightarrow U\times M(f)$.
Finally, the pullback claim will follow upon proving the universal property for the square

\[
\begin{tikzcd}
U\times M(f^s)\arrow[d,"U\times s^\lhd"']\arrow[r,"\eta^s"] & \Delta^q\arrow[d,"s"] \\
U\times M(f)\arrow[r,"\eta"'] & \Delta^p
\end{tikzcd}.
\]

Consider maps $\alpha:T\rightarrow U\times M(f),\beta:T\rightarrow\Delta^q$ fitting into a commutative diagram

\[
\begin{tikzcd}[column sep=2.5em, row sep=2.5em]
T
  \arrow[rrd, bend left=15,"\beta"]
  \arrow[ddr, bend right=15,"\alpha"']
& & \\
& U\times M(f^s) \arrow[r,"\eta^s"] \arrow[d,"U\times s^\lhd"'] & \Delta^q \arrow[d,"s"] \\
& U\times M(f) \arrow[r,"\eta"'] & \Delta^p
\end{tikzcd}.
\]

We define the map $\psi:T\rightarrow U\times M(f^s)$ by

\[
\psi(x)=
\lf\{
\begin{array}{ll}
(v,\rho(\beta(x),\eta(u,t,1)),r) & \alpha(x)=(v,t,r),r>0 \\
(v,\rho(\beta(x),\eta(u,y))) & \alpha(x)=(v,y),y\in Y.
\end{array}
\ri.
\]

The continuity of $\psi$ comes from the continuity of $\eta$ and $\rho$. Moreover, for $x\in T$ such that $\alpha(x)=(v,t,r),r>0$, we have
\begin{align*}
\eta^s\circ\psi(x) &= \eta^s\lf(v,\rho\lf(\beta(x),\eta\lf(u,t,1)\ri)\ri),r\ri) \\
&= \rho\lf(\rho\lf(\beta(x),\eta\lf(u,t,1\ri)\ri),\eta(v,t,r)\ri) \\
&= \rho(\beta(x),\eta(\alpha(x))) \\
&= \rho(\beta(x),s(\beta(x))) \\
&= \beta(x).
\end{align*}

We also have

\begin{align*}
(U\times s^\lhd)(\psi(x)) &= (U\times s^\lhd)\lf(v,\rho\lf(\beta(x),\eta\lf(u,t,1\ri)\ri),r\ri) \\
&= (v,t,r) \\
&= \alpha(x).
\end{align*}

The case $\alpha(x)=(v,y),y\in Y$ is done similarly.

Moreover, such an extension is unique because $\eta^s$ is injective.
Thus it is the unique map making the whole diagram commute. Therefore the square is a pullback. Finally, since $\eta:U\times M(f)\rightarrow\Delta^p$ is an open embedding and the pullback of an open embedding is again an open embedding, $\eta^s:U\times M(f^s)\rightarrow\Delta^q$ is an open embedding as well.

\end{proof}

\begin{lemma}\label{compatibility_faces}
Consider a commutative square
\[
\begin{tikzcd}
\Delta^{q'}\arrow[d,"s'"']\arrow[r,"d'"] & \Delta^q\arrow[d,"s"] \\
\Delta^{p'}\arrow[r,"d"'] & \Delta^p
\end{tikzcd}
\]
where horizontal maps are face maps and vertical maps are degeneracies. Consider commutative squares
\[
\begin{tikzcd}
W'\arrow[r]\arrow[d,"f'"'] & W\arrow[d,"f"] \\
Y'\arrow[r] & Y
\end{tikzcd},
\]
\[
\begin{tikzcd}
U\times M(f')\arrow[d,"\eta'"']\arrow[r] & U\times M(f) \arrow[d,"\eta"] \\
\Delta^{p'}\arrow[r,"d"'] & \Delta^p 
\end{tikzcd}
\]
and let $((f')^{s'}:(W')^{s'}\rightarrow (Y')^{s'}),(f^s:W^s\rightarrow Y^s)$ be the liftings obtained from \Cref{lifting_degeneracies}. Then there exists a natural map $\delta:((f')^{s'}:(W')^{s'}\rightarrow (Y')^{s'})\rightarrow (f^s:W^s\rightarrow Y^s)$ providing a commutative square
\[
\begin{tikzcd}
U\times M((f')^{s'})\arrow[d,"(\eta')^{s'}"']\arrow[r,"U\times\delta^\lhd"] & U\times M(f^s) \arrow[d,"\eta^s"] \\
\Delta^{q'}\arrow[r,"d'"'] & \Delta^q 
\end{tikzcd}.
\]
Moreover, such maps are compatible with compositions of face maps.
\end{lemma}

\begin{proof}
The spaces $(W')^{s'}$ and $W^s$ are defined by $(s')^{-1}(\{u\}\times W'\times\{1\})$ and $s^{-1}(\{u\}\times W\times\{1\})$, respectively, for a fixed basepoint $u\in U$. Because $s\circ d'=d\circ s'$, the image $d'((W')^{s'})$ is contained in $W^s$, so that the face map $d':\Delta^{q'}\hookrightarrow\Delta^q$ induces the map $\delta^W:(W')^{s'}\rightarrow W^s$ by restriction and similarly for $\delta^Y$. 

Now let $\rho_s,\rho_{s'}$ the maps from \Cref{translating_fibers} induced respectively by $s,s'$. From the coordinate description of these maps and the identity $s\circ d'=d\circ s'$, one deduces the equality
\[
d'(\rho_{s'}(t,x))=\rho_s(d't,dx).
\]

Because of this, the map $(f')^{s'}:(W')^{s'}\rightarrow (Y')^{s'}$ is the restriction of the map $f^s:W^s\rightarrow Y^s$ along the face $d':\Delta^{q'}\hookrightarrow\Delta^q$. This gives the commutativity of the square

\[
\begin{tikzcd}
(W')^{s'}\arrow[r,"\delta^W"]\arrow[d,"(f')^{s'}"'] & W^s\arrow[d,"f^s"] \\
(Y')^{s'}\arrow[r,"\delta^Y"'] & Y^s
\end{tikzcd}.
\]

Moreover, since the left map is simply given by restriction to a face, this also implies compatibility with compositions of face maps.

Now consider the cubical diagram:
\[
\begin{tikzcd}
U\times M((f')^{s'})\arrow[dd,"s'"',near start]\arrow[rr,dotted]\arrow[dr,"(\eta')^{s'}"'] & & U\times M(f^s)\arrow[dd,"s",near start]\arrow[dr,"\eta^s"] & \\
& \Delta^{q'}\arrow[dd,crossing over,"s'"',near start]\arrow[rr,crossing over,"d'",near start] & & \Delta^q\arrow[dd,"s",near start] \\
U\times M(f')\arrow[dr,"\eta'"']\arrow[rr] & & U\times M(f)\arrow[dr,"\eta"] & \\
& \Delta^{p'}\arrow[rr,"d"',near start] & & \Delta^p  
\end{tikzcd}.
\]
By \Cref{lifting_degeneracies}, the left and right faces are pullback squares, and thus the dotted arrow $U\times\delta^\lhd$ is induced by the universal property. Because this map is induced by the universal property, compatibility with compositions of face maps follows.

\end{proof}

\subsection{The inductive step}

\begin{lemma}\label{embedding_boundary}
Let $W$ be a simplicial set and $f:Z\rightarrow\norm{W}$ a map from a topological space $Z$. For any simplex $\sigma\in W$, we have a homeomorphism
\[
\colim{\tau\in\partial\sigma}{\Delta^\tau\times_{\norm{W}}Z}\xrightarrow{\cong}\partial\Delta^\sigma\times_{\norm{W}}Z.
\]
\end{lemma}

\begin{proof}
For each face $\tau\in\partial\sigma$, the face inclusion $\Delta^\tau\rightarrow \partial\Delta^\sigma$ restricts to a map 
\[
\psi_\tau:\Delta^\tau\times_{\norm{W}}Z\rightarrow\partial\Delta^\sigma\times_{\norm{W}}Z.
\]
This is compatible with all face inclusions, hence induces a continuous map 
\[
\psi:\colim{\tau\in\partial\sigma}{\Delta^\tau\times_{\norm{W}}Z}\rightarrow\partial\Delta^\sigma\times_{\norm{W}}Z.
\]
This map is surjective: given a point $x\in\partial\Delta^\sigma\times_{\norm{W}}Z$, it is contained in some face $\tau\subset\partial\Delta^\sigma$. Thus it is in the image of the restriction $\Delta^\tau\times_{\norm{W}} Z\rightarrow\partial\Delta^\sigma\times_{\norm{W}}Z$. Now let $x,x'$ be points such that $\psi(x)=\psi(x')=(a,b)\in\partial\Delta^\sigma\times_{\norm{W}}Z$. Fix representatives $(w,z)$ and $(w',z')$ in $\Delta^\tau\times_{\norm{W}}Z$ and $\Delta^{\tau'}\times_{\norm{W}}Z$, respectively, for some faces $\tau,\tau'\in\partial\sigma$. Necessarily $z=z'$; moreover, the images of the points $w,w'$ in $\partial\Delta^\sigma$ must be in a common face $\tau''$ of $\tau,\tau'$. This identifies the points $(w,z),(w',z')$ in the colimit, which proves injectivity. Thus $\psi$ is bijective.

Moreover, $\psi$ is a closed map. Let $V\subset\colim{\tau\in\partial\sigma}{\Delta^\tau\times_{\norm{W}}Z}$ be a closed subset. For each $\tau\in\partial\sigma$, let $V_\tau\subset\Delta^\tau\times_{\norm{W}}Z$ be the corresponding preimages. Because all inclusions 
\[
j_\tau:\Delta^\tau\times_{\norm{W}}Z\hookrightarrow\partial\Delta^\sigma\times_{\norm{W}}Z
\]
 form a closed cover of $\partial\Delta^\sigma\times_{\norm{W}}Z$, we get that each subset $j_\tau(V_\tau)\subset\partial\Delta^\sigma\times_{\norm{W}}Z$ is closed. Thus $\psi(V)=\cup_{\tau\in\partial\sigma}j_\tau(V_\tau)$ is closed since it is a finite union of closed subsets.
\end{proof}

We can now go back to the main proof. Recall the locally compact open subset $U\subset \overset{\circ}{\Delta}{}^{\sigma_k}$ containing $x$.

For $j=k$, define
\begin{enumerate}
\item[(IH1)] $(L^k\rightarrow\ast):=(\emptyset\rightarrow\ast)$ and $\iota^k:U\times C(\emptyset)\cong U\hookrightarrow\overset{\circ}{\Delta}{}^{\sigma_k}\hookrightarrow \norm{X^k}$,
\item[(IH2)] the functor $f^k:\Delta^{inj}(X^k)\rightarrow\text{Arr}(\Top)$ by $(f^k_{\sigma_k}=\emptyset\rightarrow\ast),f^k_{s\sigma_k}$ defined by \Cref{lifting_degeneracies} for any degeneracy $s\sigma_k$ of $\sigma_k$, and $(f^k_\tau=\emptyset\rightarrow\emptyset)$ otherwise,
\item[(IH3)] the cocone over $(\emptyset\rightarrow\ast)$ defined trivially.
\item[(IH4)] the natural transformation defined trivially for $\sigma_k$, defined by \Cref{lifting_degeneracies} for $s\sigma_k$ a degeneracy of $\sigma_k$ and trivially otherwise. 
\end{enumerate}

Now assume the construction holds at step $j-1$. Using (IH2), define 
\[
(f^{\partial}:W^{\partial}\rightarrow Y^{\partial}):=\colim{\partial\sigma_j}{f^{j-1}}.
\]

By (IH3), we get an induced map $\gamma^\partial:(f^\partial:W^\partial\rightarrow Y^\partial)\rightarrow (L^{j-1}\rightarrow\ast)$.

In addition, for each proper face $\tau$ of $\sigma_j$, we have by (IH4) a pullback square

\[
\begin{tikzcd}
U\times M(f^{j-1}_\tau)\arrow[d,"\eta^{j-1}_\tau"']\arrow[r,"(\gamma^{j-1}_\tau)^\lhd"]\arrow[dr,"\lrcorner",phantom,very near start] & U\times C(L^{j-1})\arrow[d,"\iota^{j-1}"] \\
\Delta^\tau\arrow[r] & \norm{X^{j-1}}
\end{tikzcd}
\]

and the map $\iota^{j-1}:U\times C(L^{j-1})\rightarrow \norm{X^{j-1}}$ is an open embedding by (IH1). Thus the map $\eta^{j-1}_\tau:U\times M(f^{j-1}_\tau)\rightarrow\Delta^\tau$ is also an open embedding. Because mapping cylinders commute with colimits, we then get 

\begin{align*}
U\times M(f^\partial) & \cong U\times M(\colim{\partial\sigma_j}{f^{j-1}}) \\
& \cong\colim{\tau\in\partial\sigma_j}{U\times M(f^{j-1}_\tau)} \\
& \cong\colim{\tau\in\partial\sigma_j}{\Delta^\tau\times_{\norm{X^{j-1}}}U\times C(L^{j-1})} \\
& \cong\partial\Delta^{\sigma_j}\times_{\norm{X^{j-1}}}U\times C(L^{j-1}).
\end{align*}

This gives an embedding $\eta^\partial:U\times M(f^\partial)\hookrightarrow\partial\Delta^{\sigma_j}$ by \Cref{embedding_boundary}. Moreover, since each $W_\tau,Y_\tau$ is compact by (IH2) and the colimit is taken over a finite category, we deduce that $W^\partial,Y^\partial$ are compact as well.

 We then define $(f^j_{\sigma_j}:W_{\sigma_j}\rightarrow Y_{\sigma_j}):=(\bar{M}(f^\partial:W^\partial\rightarrow Y^\partial)\rightarrow Y^\partial)$.

In addition, \Cref{prism_lemma} provides a pushout square (note that we use the fact that $U$ is locally compact)

\[
\begin{tikzcd}
U\times M(f^\partial)\arrow[r]\arrow[d,hook] & U\times C(L^{j-1})\arrow[d,"U\times (l^j)^\lhd"] \\
U\times M(r_{f^\partial})\arrow[r,"U\times (\gamma^j_{\sigma_j})^\lhd"'] & U\times C(L^j)
\end{tikzcd}
\]

where $L^j$ fits into a pushout square

\[
\begin{tikzcd}
W^\partial\arrow[r,"\gamma^\partial"]\arrow[d,"W^\partial\times 1"'] & L^{j-1}\arrow[d,"l^j"] \\
\bar{M}(f^\partial)\arrow[r,"\gamma^j_{\sigma_j}"'] & L^j\arrow[ul,"\ulcorner",very near start,phantom]
\end{tikzcd}.
\]

Since $W^\partial,Y^\partial,L^{j-1}$ are compact, we deduce that $L^j$ is compact as well.

By \Cref{cone_prism}, the inclusion $U\times M(f^\partial)\rightarrow U\times M(r_{f^\partial})$ is identified with the inclusion $U\times M(f^\partial)\times\{0\}\hookrightarrow U\times M(f^\partial)\times\R_{\geq 0}$. This fits into a commutative square
\[
\begin{tikzcd}
U\times M(f^\partial)\times\{0\}\arrow[d,hook]\arrow[r,hook,"\iota^\partial"] & \partial\Delta^{\sigma_j}\arrow[d,hook] \\
U\times M(f^\partial)\times\R_{\geq 0}\arrow[r,hook] & \partial\Delta^{\sigma_j}\times\R_{\geq 0}
\end{tikzcd}
\]
and composing with a collar embedding $\partial\Delta^{\sigma_j}\times\R_{\geq 0}\hookrightarrow\Delta^{\sigma_j}$ defines an embedding $U\times M(f^\partial)\times\R_{\geq 0}\hookrightarrow\Delta^{\sigma_j}$ which is open since $\eta^\partial:U\times M(f^\partial)\hookrightarrow\partial\Delta^{\sigma_j}$ is an open embedding.

Thus the pushout defining $U\times C(L^j)$ maps into another pushout square

\[
\begin{tikzcd}
\partial\Delta^{\sigma_j}\arrow[d]\arrow[r] & \norm{X^{j-1}}\arrow[d] \\
\Delta^{\sigma_j}\arrow[r] & \norm{X^j}
\end{tikzcd}.
\]

This defines the map $\iota^j:U\times C(L^j)\rightarrow \norm{X^j}$. 

We now prove that the induction hypotheses (IH1), (IH2), (IH3), and (IH4) still hold.

\begin{lemma}\label{injective_pushout}
Consider a cubical diagram of topological spaces
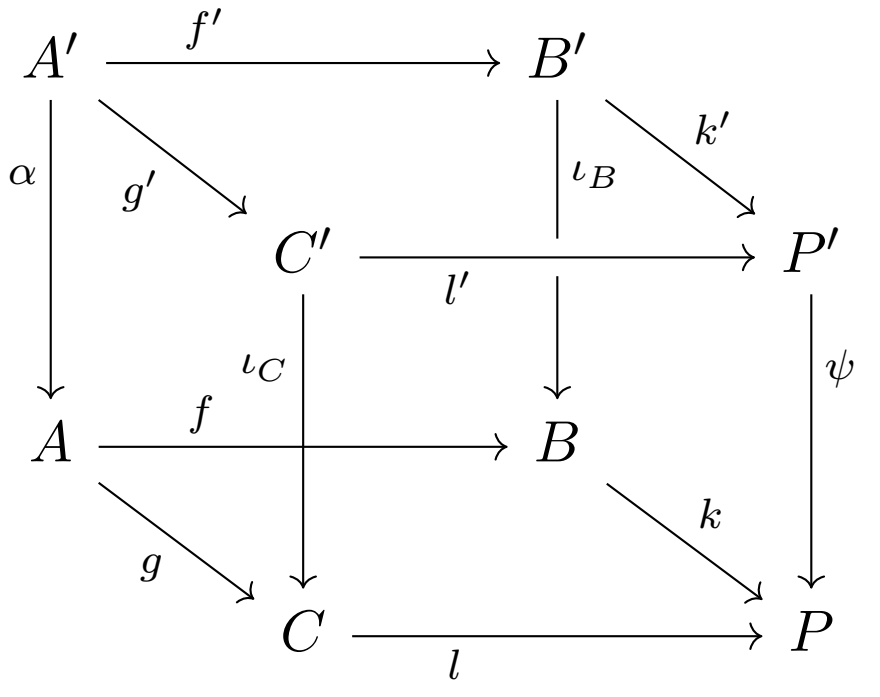
\[
\begin{tikzcd}
A'\arrow[dd,"\alpha"',near start]\arrow[rr,"f'",near start]\arrow[dr,"g'"'] & & B'\arrow[dd,"\iota_B",near start]\arrow[dr,"k'"] & \\
& C'\arrow[dd,"\iota_C"',near start,crossing over]\arrow[rr,"l'"',near start,crossing over] & & P'\arrow[dd,"\psi",near start]  \\
A\arrow[dr,"g"']\arrow[rr,"f",near start] & & B\arrow[dr,"k"] & \\
& C\arrow[rr,"l"',near start] & & P  
\end{tikzcd}
\]

where the top and bottom faces are pushouts and the left and back faces are pullbacks. Assume that the maps $\iota_B,\iota_C$ are open embeddings. Then the induced map $\psi:P'\rightarrow P$ is an open embedding. Moreover, the front and right faces are pullbacks.
\end{lemma}

\begin{proof}
Let $a'\in A',a\in A$ such that $f(\alpha(a'))=f(a)$. In particular $f(a)=\iota_B(f'(a'))$. Since the back square is a pullback, this provides an element $\bar{a}\in A'$ such that $\alpha(\bar{a})=a$. The same argument shows that $g(\alpha(a'))=g(a)$ implies $a\in\alpha(A')$.

Moreover, a pair $(a,b')\in A\times B'$ such that $f(a)=\iota_B(b')$ provides an element $a'\in A'$ such that $\alpha(a')=a,f'(a')=b'$ and likewise for $C$. Writing $\sim_A,\sim_{A'}$ the equivalence relations defining respectively the pushouts $(B\sqcup C)/\sim_A\cong P, (B'\sqcup C')/\sim_{A'}$, both observations imply that for any $x_0,x_1\in B'\sqcup C'$, we have $(x_0\sim_{A'}x_1)$ whenever $(\iota_B\sqcup \iota_C)(x_0)\sim_A(\iota_B\sqcup \iota_C)(x_1))$ and $x_0\sim_A(\iota_B\sqcup \iota_C)(x_1)$ only if $x_1\in(\iota_B\sqcup \iota_C)(B'\sqcup C')$. We then deduce that the induced map $\psi:P'\rightarrow P$ is injective.

Next we prove that the front face is a pullback. Since $\iota_C$ is an embedding, hence defines a homeomorphism $C'\xrightarrow{\cong} \iota_C(C')$ and $\psi$ is injective, the claim will follow upon proving that given $c\in C,p\in P'$ such that $l(c)=\psi(p)$, then $c\in\iota_C(C')$. Let $x\in B'\sqcup C'$ identified to $p$ via the quotient $(B'\sqcup C')/\sim_{A'}\cong P'$. Then $c\sim_A \iota_{B\sqcup C}(x)$. By the previous observation, this implies $c\in\iota_{B\sqcup C}(B'\sqcup C')$ and thus $c\in\iota_C(C')$ as desired. The same argument shows the right face is a pullback.

Finally, because $\psi$ is injective, to prove that $\psi$ is an open embedding, it is enough to prove it is an open map. Let $U\subset P'$ be an open subset. Since the front face is a pullback, we have $l^{-1}(\psi(U))=\iota_C((l')^{-1}(U))$. The latter is open since $l'$ is continuous and $\iota_C$ is an open map. Likewise, we deduce that $k^{-1}(\psi(U))$ is open in $B$. By definition of the quotient topology defining the pushout, we deduce that $\psi(U)$ is open in $P$.
\end{proof}

\begin{lemma}\label{some_pullback_squares}
The following two commutative squares are pullbacks:
\[
\begin{tikzcd}
U\times M(f^\partial)\arrow[d]\arrow[r] & \partial\Delta^{\sigma_j} \arrow[d] \\
U\times M(f^j_{\sigma_j})\arrow[r] & \Delta^{\sigma_j}
\end{tikzcd},
\begin{tikzcd}
U\times M(f^\partial)\arrow[d,hook,"\eta^\partial"']\arrow[r,"\gamma^\partial"] & U\times C(L^{j-1})\arrow[d,hook,"\iota^{j-1}"] \\
\partial\Delta^{\sigma_j}\arrow[r,"q^\partial"'] & \norm{X^{j-1}}
\end{tikzcd}.
\]
\end{lemma}

\begin{proof}
By definition of collar embedding defining the map $U\times M(f^j_{\sigma_j})\hookrightarrow\Delta^{\sigma_j}$ , we have a first pullback square

\[
\begin{tikzcd}
U\times M(f^\partial)\arrow[r]\arrow[d]\arrow[dr,"\lrcorner",phantom,very near start] & \partial\Delta^{\sigma_j}\arrow[d] \\
U\times M(f^j_{\sigma_j})\arrow[r] & \partial\Delta^{\sigma_j}\times \R_{\geq 0}
\end{tikzcd}.
\]

We also have another pullback square

\[
\begin{tikzcd}
\partial\Delta^{\sigma_j}\arrow[r,equal]\arrow[d]\arrow[dr,"\lrcorner",phantom,very near start] & \partial\Delta^{\sigma_j}\arrow[d] \\
\partial\Delta^{\sigma_j}\times \R_{\geq 0}\arrow[r] & \Delta^{\sigma_j}
\end{tikzcd}.
\]

Composing both gives the first pullback square. 

Next we show that the square
\[
\begin{tikzcd}
U\times M(f^\partial)\arrow[d,hook,"\eta^\partial"']\arrow[r,"\gamma^\partial"] & U\times C(L^{j-1})\arrow[d,hook,"\iota^{j-1}"] \\
\partial\Delta^{\sigma_j}\arrow[r,"q^\partial"'] & \norm{X^{j-1}}
\end{tikzcd}
\]
is a pullback. Since both vertical maps are embeddings, it is enough to prove that for any element $x\in\partial\Delta^{\sigma_j},q^{\partial}(x)\in\iota^{j-1}(U\times C(L^{j-1}))$ if and only if $x\in \eta^{\partial}(U\times M(f^\partial))$. The "if" direction follows from the commutativity of the square. Conversely, fix $x\in\partial\Delta^{\sigma_j}$ and assume that $q^\partial(x)\in\iota^{j-1}(U\times C(L^{j-1}))$. Let $\tau$ be a strict face of $\sigma_j$ such that the image of the face inclusion $\Delta^\tau\hookrightarrow\partial\Delta^{\sigma_j}$ contains $x$. We get a commutative rectangle
\[
\begin{tikzcd}
U\times M(f^{j-1}_\tau)\arrow[d,"\eta^{j-1}_\tau"']\arrow[r] & U\times M(f^\partial)\arrow[r,"\gamma^\partial"]\arrow[d,"\eta^\partial"'] & U\times C(L^{j-1})\arrow[d,"\iota^{j-1}"] \\
\Delta^\tau\arrow[r,hook] & \partial\Delta^{\sigma_j}\arrow[r,"q^\partial"'] & \norm{X^{j-1}}
\end{tikzcd}.
\]
Since the whole rectangle is a pullback square by (IH4) and the vertical maps are embeddings, by identifying $x$ with a point $x'$ in $\Delta^\tau$, we deduce that $x'\in\eta^{j-1}_\tau(U\times M(f^{j-1}_\tau))$. The commutativity of the left square further implies that $x\in\eta^\partial(U\times M(f^\partial))$ as desired. 
\end{proof}

\begin{lemma}\label{pb_open_embedding}
The following commutative squares are pullbacks:
\[
\begin{tikzcd}
U\times M(f^j_{\sigma_j})\arrow[d,hook]\arrow[r] & U\times C(L^j) \arrow[d,"\iota^j"] \\
\Delta^{\sigma_j}\arrow[r] & \norm{X^j}
\end{tikzcd},
\begin{tikzcd}
U\times C(L^{j-1})\arrow[d,hook,"\iota^{j-1}"']\arrow[r] & U\times C(L^j) \arrow[d,"\iota^j"] \\
\norm{X^{j-1}}\arrow[r] & \norm{X^j}
\end{tikzcd}.
\]
Moreover the map $\iota^j:U\times C(L^j)\rightarrow \norm{X^j}$ is an open embedding.
\end{lemma}

\begin{proof}

Consider the cube
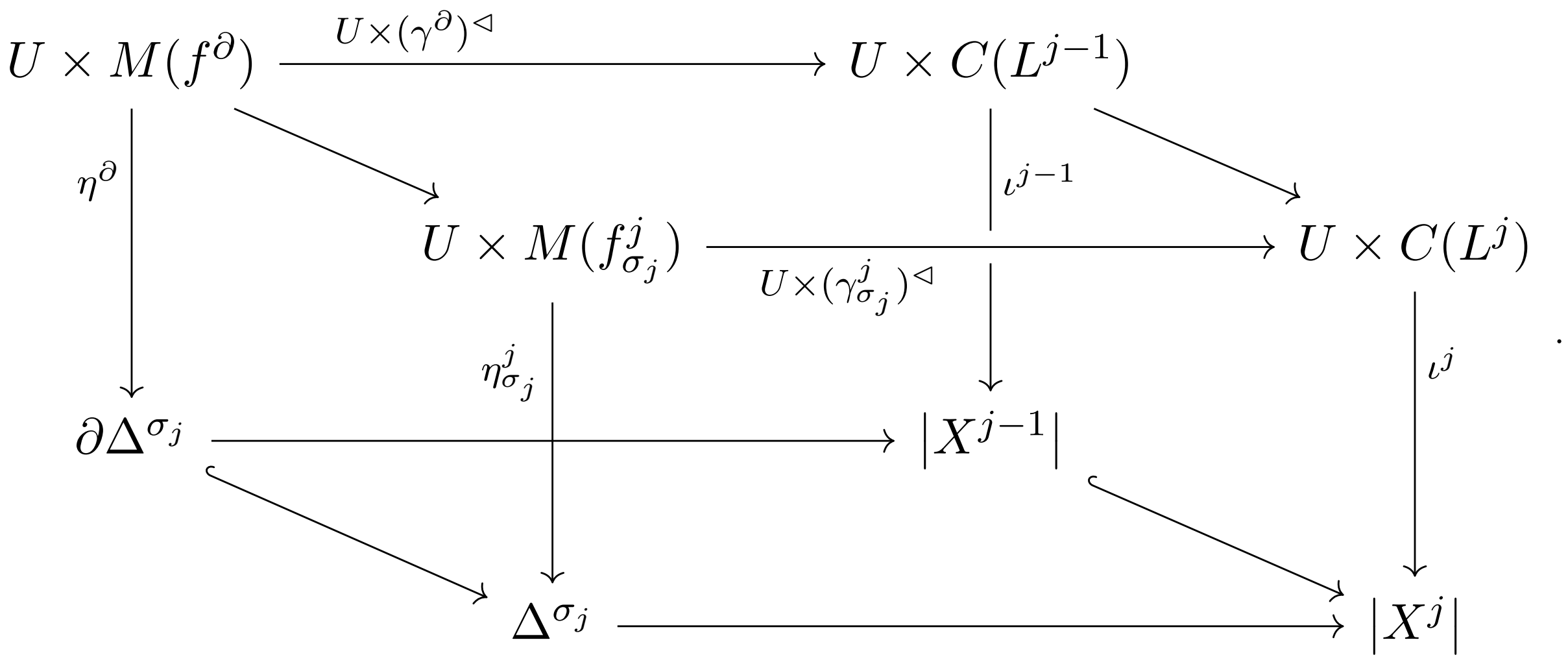
\[
\resizebox{\textwidth}{!}{%
\begin{tikzcd}[ampersand replacement=\&]
U\times M(f^\partial)
  \arrow[dr]
  \arrow[dd,near start,"\eta^\partial"']
  \arrow[rr,near start,"U\times (\gamma^\partial)^\lhd"]
\& \&
U\times C(L^{j-1})
  \arrow[dr]
  \arrow[dd,near start,"\iota^{j-1}"]
\&
\\
\&
U\times M(f^j_{\sigma_j})
  \arrow[dd,near start,"\eta^j_{\sigma_j}"',crossing over]
  \arrow[rr,near start,"U\times (\gamma^j_{\sigma_j})^\lhd"',crossing over]
\& \&
U\times C(L^j)
  \arrow[dd,"\iota^j",near start]
\\
\partial\Delta^{\sigma_j}
  \arrow[dr,hook]
  \arrow[rr]
\& \&
\norm{X^{j-1}}
  \arrow[dr,hook]
\&
\\
\&
\Delta^{\sigma_j}
  \arrow[rr]
\& \&
\norm{X^j}
\end{tikzcd}%
}.
\]
The top and bottom faces are pushouts by construction. The left and back faces are pullbacks by \Cref{some_pullback_squares}. The maps $\eta^j_{\sigma_j}:U\times M(f^j_{\sigma_j})\hookrightarrow\Delta^{\sigma_j},\iota^{j-1}:U\times C(L^{j-1})\hookrightarrow \norm{X^{j-1}}$ are open embeddings. By \Cref{injective_pushout}, the claim follows.
\end{proof}

\begin{lemma}\label{IH1_stratification}
The topological space $L^j$ comes with a stratification $\pi^j:L^j\rightarrow\pos{X^j}_{\sigma_k<}$ such that the map $\iota^j:U\times C(L^j)\rightarrow \norm{X^j}$ fits into a commutative square
\[
\begin{tikzcd}
U\times C(L^j)\arrow[r,"\iota^j"]\arrow[d,"(\pi^j)^{\lhd}"'] & \norm{X^j}\arrow[d,"\pi_{X^j}"] \\
\pos{X^j}_{\sigma_k\leq}\arrow[r] & \pos{X^j}
\end{tikzcd}.
\]
\end{lemma}

\begin{proof}
By construction, $L^j$ fits into a pushout square

\[
\begin{tikzcd}
W^\partial\arrow[r,"\gamma^\partial"]\arrow[d,"W^\partial\times 1"'] & L^{j-1}\arrow[d,"l^j"] \\
\bar{M}(f^\partial)\arrow[r,"(\gamma^j_{\sigma_j})^\lhd"'] & L^j\arrow[ul,"\ulcorner",very near start,phantom]
\end{tikzcd}.
\]

By (IH1), we have the stratification $\pi^{j-1}:L^{j-1}\rightarrow\pos{X^{j-1}}_{\sigma_k<}\hookrightarrow\pos{X^j}_{\sigma_k<}$. Moreover we have another pushout 

\[
\begin{tikzcd}
W^\partial\arrow[r,"f^\partial"]\arrow[d,"W^\partial\times 0"'] & Y^\partial\arrow[d] \\
W^\partial\times [0,1]\arrow[r] & \bar{M}(f^\partial)\arrow[ul,"\ulcorner",very near start,phantom]
\end{tikzcd}.
\]

We define a stratification $\kappa:W^\partial\times [0,1]\rightarrow\pos{X^j}_{\sigma_k<}$ by

\[
\kappa(w,r) =
\lf\{
\begin{array}{ll}
\sigma_j & r<1, \\
\lambda^{j-1}(\gamma^\partial(w)) & r=1. 
\end{array}
\ri.
\]

This map is indeed continuous: by construction of $W^\partial$, for any $w\in W^\partial$ we have $\lambda^{j-1}(\gamma^\partial(w))\leq\sigma_j$. Thus, given an open subset $P\subset\pos{X^j}_{\sigma_k<}$, we have
\[
\kappa^{-1}(P) =
\lf\{
\begin{array}{ll}
W^\partial\times [0,1[\cup \lf((\lambda^{j-1}\circ \gamma^\partial)^{-1}(P)\times\{1\}\ri) & \sigma_j\in P, \\
\emptyset & \text{otherwise}
\end{array}
\ri.
\]
This is open in both cases. Then, the pair $(\kappa,\cst_{\sigma_j})$ defines the stratification $\bar{\pi}:\bar{M}(f^\partial)\rightarrow\pos{X^j}_{\sigma_k<}$, where $\cst_{\sigma_j}:Y^\partial\rightarrow\pos{X^j}$ is the constant stratification sending all points to $\sigma_j$. Because both stratifications $\bar{\pi}:\bar{M}(f^\partial)\rightarrow\pos{X^j}_{\sigma_k<}$ and $\pi^{j-1}:L^{j-1}\rightarrow\pos{X^j}_{\sigma_k<}$ satisfy $\pi^{j-1}\circ \gamma^{\partial}=\bar{\pi}\circ (W^\partial\times 1)$, we get a stratification $\pi^j:L^j\rightarrow\pos{X^j}_{\sigma_k<}$.

It remains to check the compatibility of this stratification with the map $\iota^j:U\times C(L^j)\rightarrow \norm{X^j}$. Since the stratification $\norm{X^j}\rightarrow\pos{X^j}$ is the universal map provided by the pushout to the commutative square

\[
\begin{tikzcd}
\partial\Delta^{\sigma_j}\arrow[r]\arrow[d] & \norm{X^{j-1}}\arrow[d,"\pi_{X^{j-1}}"] \\
\Delta^{\sigma_j}\arrow[r] & \pos{X^j}
\end{tikzcd},
\]

and the map $(\pi^j)^\lhd:U\times C(L^j)\rightarrow \pos{X^j}_{\sigma_k\leq }$ is itself induced by universality applied to
\[
\begin{tikzcd}
U\times M(f^\partial)\arrow[d,"U\times (W^\partial\times 1)^\lhd"']\arrow[r,"U\times (\gamma^\partial)^\lhd"] & U\times C(L^{j-1})\arrow[d,"(\pi^{j-1})^\lhd"] \\
U\times M(f^j_{\sigma_j})\arrow[r,"U\times\bar{\pi}^\lhd"'] & \pos{X^j}_{\sigma_k\leq}
\end{tikzcd},
\]

the compatibility will follow upon proving that the following cube is commutative:

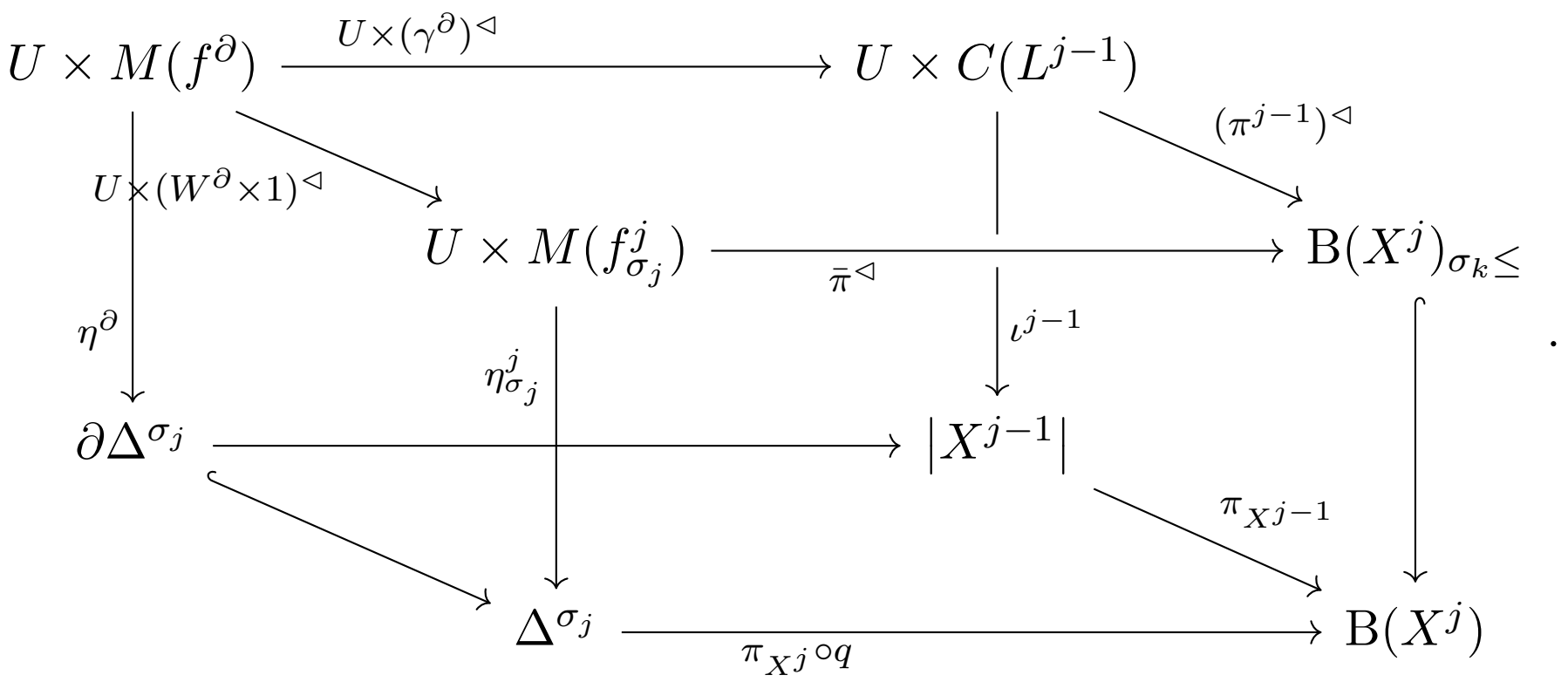
\[
\begin{tikzcd}
U\times M(f^\partial)\arrow[dr,"U\times (W^\partial\times 1)^\lhd"']\arrow[dd,"\eta^\partial"',near end]\arrow[rr,"U\times(\gamma^\partial)^\lhd",near start] & & U\times C(L^{j-1})\arrow[dd,"\iota^{j-1}",near end]\arrow[dr,"(\pi^{j-1})^\lhd"] & \\
& U\times M(f^j_{\sigma_j})\arrow[dd,"\eta^j_{\sigma_j}"',near start,crossing over]\arrow[rr,"\bar{\pi}^\lhd"',near start,crossing over] & & \pos{X^j}_{\sigma_k\leq}\arrow[dd,hook] \\
\partial\Delta^{\sigma_j}\arrow[dr,hook]\arrow[rr] & & \norm{X^{j-1}}\arrow[dr,"\pi_{X^{j-1}}"] & \\
& \Delta^{\sigma_j}\arrow[rr,"\pi_{X^j}\circ q"',near start] & & \pos{X^j}
\end{tikzcd}.
\]

The top and bottom faces commute by construction, the right face commutes by (IH1), and the back face commutes by (IH4). By construction, $\eta^j_{\sigma_j}:U\times M(f^j_{\sigma_j})\rightarrow\Delta^{\sigma_j}$ is identified with a collar embedding inside $\Delta^{\sigma_j}$ of the embedding $\eta^\partial:U\times M(f^\partial)
\hookrightarrow\partial\Delta^{\sigma_j}$, using the homeomorphism 
\begin{align*}
& h:\bigl(U\times (W^\partial\times 1)^{\lhd}:U\times M(f^\partial)\hookrightarrow U\times M(f^j_{\sigma_j})\bigr) \\
& \cong \bigl((U\times M(f^\partial))\times\{0\}\hookrightarrow (U\times M(f^\partial))\times\R_{\geq 0}\bigr)
\end{align*}
from \Cref{cone_prism}; in particular, the left square commutes. Moreover, \Cref{cone_prism} provides a commutative diagram
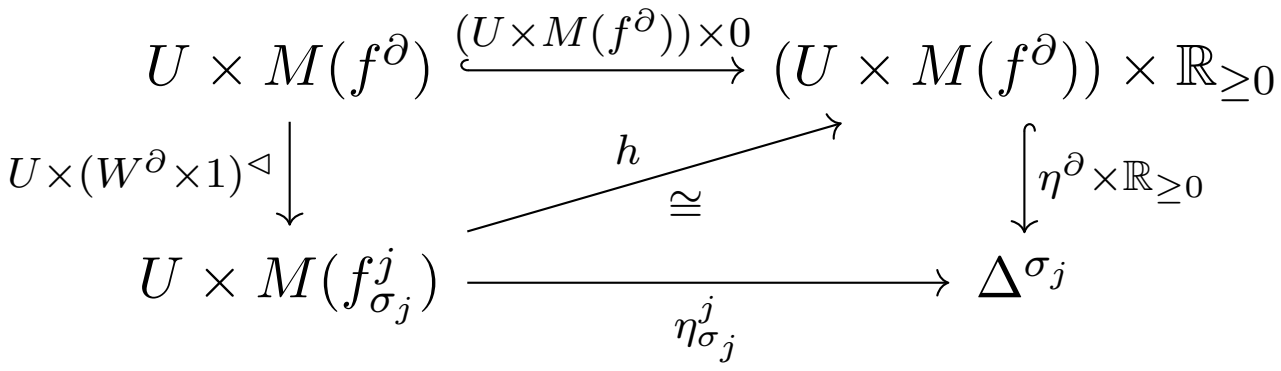
\[
\begin{tikzcd}
U\times M(f^\partial)\arrow[d,"U\times (W^\partial\times 1)^\lhd"']\arrow[rr,hook,"(U\times M(f^\partial))\times 0"] & & (U\times M(f^\partial))\times\R_{\geq 0}\arrow[d,"\eta^\partial\times\R_{\geq 0}",hook] \\
U\times M(f^j_{\sigma_j})\arrow[rr,"\eta^j_{\sigma_j}"']\arrow[urr,"h","\cong"'] & & \Delta^{\sigma_j}
\end{tikzcd}
\]
where $\eta^\partial\times\R_{\geq 0}$ denotes a collar embedding of the embedding $\eta^\partial$ inside $\Delta^{\sigma_j}$. By definition of collar embedding and the stratification $(\pi_{X^j}\circ q):\Delta^{\sigma_j}\rightarrow\pos{X^j}$, we have 
\[
(\pi_{X^j}\circ q\circ (\eta^\partial\times\R_{\geq 0}))(x)=\sigma_j
\]
for any $x\in (U\times M(f^\partial))\times\R_{>0}$ and 
\[
(\pi_{X^j}\circ q\circ (\eta^\partial\times\R_{\geq 0}))(x)\neq\sigma_j
\] 
for any $x\in (U\times M(f^\partial))\times\{0\}$. This implies that for $x\in U\times M(f^j_{\sigma_j})$, 
\[
(\pi_{X^j}\circ q\circ \eta^j_{\sigma_j})(x)=\sigma_j
\]
 if and only if 
\[ 
 x\not\in (U\times (W^\partial\times 1)^\lhd)(U\times M(f^\partial)).
 \]
 Moreover, if $x=(U\times (W^\partial\times 1)^\lhd)(w)$ for some $w\in U\times M(f^\partial)$, then 
 \[
 (\pi_{X^j}\circ q\circ \eta^j_{\sigma_j})(x)=(\pi_{X^j}\circ q\circ\eta^\partial)(w)=\lf((\pi^{j-1})^\lhd\circ (U\times(\gamma^\partial)^\lhd))\ri)(w).
 \]
Because of these identities, we deduce that $(\pi_{X^j}\circ q\circ\eta^j_{\sigma_j})=\bar{\pi}^\lhd$, hence the front face commutes.

\end{proof}

\begin{lemma}\label{IH2}
There is a functor $f^j:\Delta^{inj}(X^j)\rightarrow\text{Arr}(\Top)$ such that $f^j_{\sigma_j}=(r_{f^\partial}:\bar{M}(f^\partial)\rightarrow Y^\partial)$ and $f^j_{|\Delta^{inj}(X^{j-1})}=f^{j-1}$. 
\end{lemma}

\begin{proof}
Given a strict face $d:\tau\rightarrow\sigma_j$, we first define 
\[
f^j(d):(f^{j-1}_\tau:W_\tau\rightarrow Y_\tau)\rightarrow(f^j_{\sigma_j}:\bar{M}(f^\partial)\rightarrow Y^\partial)
\]
 by the composition
\[
(f^{j-1}_\tau:W_\tau\rightarrow Y_\tau)\rightarrow(f^\partial:W^\partial\rightarrow Y^\partial)\xrightarrow{W^\partial\times 1}(f^j_{\sigma_j}:\bar{M}(f^\partial)\rightarrow Y^\partial)
\]
where the first map is the natural map into the colimit $\colim{\partial\sigma_j}{f^{j-1}}\cong f^{\partial}$. Since this map is natural, this ensures that given a composition of faces $\tau'\xrightarrow{d'}\tau\xrightarrow{d}\sigma_j$, we have $f^j(d\circ d')=f^j(d)\circ f^j(d')$. In addition, \Cref{lifting_degeneracies} and \Cref{compatibility_faces} ensure the extension to degeneracies of $\sigma_j$ and hence give the desired functor $f^j:\Delta^{inj}(X^j)\rightarrow\text{Arr}(\Top)$. 
\end{proof}

\begin{lemma}\label{IH3}
The functor $f^j:\Delta^{inj}(X^j)\rightarrow\text{Arr}(\Top)$ extends to a diagram over $(L^j\rightarrow\ast)$.
\end{lemma}

\begin{proof}
Recall that $L^j$ fits into a pushout square

\[
\begin{tikzcd}
W^\partial\arrow[d,"W^\partial\times 1"']\arrow[r,"\gamma^\partial"] & L^{j-1}\arrow[d] \\
\bar{M}(f^\partial)\arrow[r] & L^j\arrow[ul,"\ulcorner",very near start,phantom]
\end{tikzcd}.
\]

For each $\tau\in\Delta^{inj}(X^{j-1})$, we get a map $(W_\tau\rightarrow Y_\tau)\rightarrow (L^j\rightarrow\ast)$ induced by the composition

\[
(W_\tau\rightarrow Y_\tau)\rightarrow (L^{j-1}\rightarrow\ast)\rightarrow (L^j\rightarrow\ast).
\]

Moreover, the map $\bar{M}(f^\partial)\rightarrow L^j$ defining the pushout induces the map
\[
(f^j_{\sigma_j}:\bar{M}(f^\partial)\rightarrow Y^\partial)\rightarrow (L^j\rightarrow\ast)
\]

and the same argument as in \Cref{IH2} shows it is compatible with face restrictions. In addition, 
\Cref{lifting_degeneracies} and \Cref{compatibility_faces} extend these natural maps to a cocone $f^j\Rightarrow (L^j\rightarrow\ast)$ defined on $\Delta^{inj}(X^j)$.
\end{proof}

\begin{lemma}\label{IH4}
For each $\tau\in\Delta^{inj}(X^j)$, there is a natural map $\eta^j_\tau:U\times M(f^j_\tau)\rightarrow\Delta^\tau$ fitting into a pullback square
\[
\begin{tikzcd}
U\times M(f^j_\tau)\arrow[d,"\eta^j_\tau"']\arrow[dr,very near start,"\lrcorner",phantom]\arrow[r,"(\gamma^j_\tau)^\lhd"] & U\times C(L^j)\arrow[d] \\
\Delta^\tau\arrow[r] & \norm{X^j}
\end{tikzcd}.
\]
\end{lemma}

\begin{proof}
For any $\tau\in\Delta^{inj}(X^{j-1})$, we have by construction a natural homeomorphism 
\[
\Delta^{\tau}\times_{\norm{X^j}}(U\times C(L^j))\cong \Delta^{\tau}\times_{\norm{X^{j-1}}}(U\times C(L^{j-1}))
\]
 and the claim follows from (IH4). Moreover, from the collar embedding $\eta^j_{\sigma_j}:U\times M(f^j_{\sigma_j})\hookrightarrow\Delta^{\sigma_j}$ fitting into the commutative square

\[
\begin{tikzcd}
U\times M(f^\partial)\arrow[r,"U\times f^j(d)^\lhd"]\arrow[d] & U\times M(f^j_{\sigma_j})\arrow[d] \\
\partial\Delta^{\sigma_j}\arrow[r] & \Delta^{\sigma_j}
\end{tikzcd}.
\] 

we deduce that given a strict face $d:\tau\rightarrow\sigma_j$, there is another commutative square
\[
\begin{tikzcd}
U\times M(f^j_\tau)\arrow[r]\arrow[d] & U\times M(f^j_{\sigma_j})\arrow[d] \\
\Delta^\tau\arrow[r] & \Delta^{\sigma_j}
\end{tikzcd},
\] 

The square

\[
\begin{tikzcd}
U\times M(f^j_{\sigma_j})\arrow[d,"\eta^j_{\sigma_j}"']\arrow[r,"U\times(\gamma^j_{\sigma_j})^\lhd"] & U\times C(L^j)\arrow[d,"\iota^j"] \\
\Delta^{\sigma_j}\arrow[r,"q_{\sigma_j}"'] & \norm{X^j}
\end{tikzcd}
\]

is a pullback square. 

Finally, for degenerate simplices in $\Delta^{inj}(X^j)$, the naturality comes from \Cref{lifting_degeneracies} and \Cref{compatibility_faces}. 
\end{proof}

Combining \Cref{IH1_stratification}, \Cref{IH2}, \Cref{IH3}, and \Cref{IH4} completes the inductive step and concludes the proof of \Cref{conicality}.

\begin{corollary}[\Cref{main_strat_A}]
Let $X$ be a locally finite simplicial set. Then the stratification $\pi_X:\norm{X}\rightarrow\pos{X}$ is conical.
\end{corollary}

\begin{proof}
Let $x\in\norm{X}$ and $\sigma$ be the stratum containing it. Let $X^{\sigma}\subset X$ be the smallest simplicial subset containing all simplices for which $\sigma$ is a face and $V:=\pi_X^{-1}(\pos{X}_{\sigma\leq})$. Since $X$ is locally finite, $X^{\sigma}$ is finite. In addition, because $\pos{X}_{\sigma\leq}$ is open in $\pos{X}$, $V$ is open in $\norm{X}$. Moreover, since the inclusion $\pos{X}_{\sigma\leq}\hookrightarrow\pos{X}$ factors through $\pos{X^\sigma}$ by construction, the inclusion $V\hookrightarrow\norm{X}$ factors through an inclusion $V\hookrightarrow\norm{X^\sigma}$. Now, by \Cref{conicality}, there exists a conical neighbourhood $U\times C(L)\hookrightarrow V\subset\norm{X^\sigma}$ around $x$. Then $U\times C(L)$ is open in $V$, and because $V$ is open in $\norm{X}$, $U\times C(L)$ is open in $\norm{X}$. This defines the desired conical neighbourhood in $\norm{X}$.

\end{proof}

\section{Descent over the category of simplices}

Let $\sSet^{\text{fin}}\subset\sSet$ be the full subcategory consisting of finite simplicial sets. By \Cref{functoriality_pi} and \Cref{main_strat_A}, we have a functor $\exit{\pos{-}}{\norm{-}}:\sSet^{\text{fin}}\rightarrow\Cat$. This fits into the diagram

\[
\begin{tikzcd}
\Delta\arrow[d,"y"']\arrow[rr,"\exit{\pos{-}}{\Delta^-}"] & & \Cat \\
\sSet^{\text{fin}}\arrow[urr,"\exit{\pos{-}}{\norm{-}}"']
\end{tikzcd}
\]

where $y:\Delta\rightarrow\sSet^{\text{fin}}$ is the Yoneda embedding.

The goal of this section is to prove that this diagram exhibits a left Kan extension.

\begin{theorem}[\Cref{main_strat_B}]
Let $X$ be a finite simplicial set. Then there is an equivalence of $\infty$-categories
\[
\colim{\tau\in\Delta(X)}{\exit{\pos{\tau}}{\Delta^\tau}}\simeq\exit{\pos{X}}{\norm{X}}.
\]
\end{theorem}

Let $\Psh(\Delta)$ be the $\infty$-category of space-valued presheaves on $\Delta$ and write $\bar{y}:\Delta\hookrightarrow\Psh(\Delta)$ for the Yoneda embedding. Viewing a finite simplicial set as a discrete presheaf, we get a fully faithful inclusion $j:\sSet^{\text{fin}}\hookrightarrow\Psh(\Delta)$ providing a factorization $j\circ y\simeq\bar{y}$.

Now consider the left Kan extension

\[
\begin{tikzcd}
\Delta\arrow[rr,"\exit{\pos{-}}{\Delta^-}"]\arrow[d,"y"'] & & \Cat \\
\Psh(\Delta)\arrow[urr,"\text{Lan}_{\bar{y}} (\exit{\pos{-}}{\Delta^-})"',dotted]
\end{tikzcd}.
\]

This left Kan extension restricts to a functor $L:\sSet^{\text{fin}}\rightarrow\Cat$ which is nothing but the left Kan extension of $\exit{\pos{-}}{\Delta^-}:\Delta\rightarrow\Cat$ along the Yoneda embedding $y:\Delta\hookrightarrow\sSet^{\text{fin}}$. Hence, \Cref{main_strat_B} can be reformulated as proving that we have a natural equivalence $L\simeq\exit{\pos{-}}{\norm{-}}$.

\begin{lemma}\label{pushout_sset}
Consider a pushout square in $\sSet^{\text{fin}}$
\[
\begin{tikzcd}
A\arrow[d,hook]\arrow[r] & X\arrow[d] \\
B\arrow[r] & P\arrow[ul,"\ulcorner",very near start,phantom]
\end{tikzcd}
\]
where the map $A\hookrightarrow B$ is a monomorphism. Under the fully faithful inclusion $j:\sSet^{\text{fin}}\hookrightarrow\Psh(\Delta)$, the square becomes a pushout in $\Psh(\Delta)$.
\end{lemma}

\begin{proof}
It is enough to prove that for each $p\geq 0$, the following square
\[
\begin{tikzcd}
A_p\arrow[d,hook]\arrow[r] & X_p\arrow[d] \\
B_p\arrow[r] & P_p
\end{tikzcd}
\]

defines a colimit in $\cS$ when viewing each object as a discrete space. Alternatively, under the embedding $\Set\hookrightarrow\sSet$ sending a set to the corresponding constant simplicial set, it is enough to prove that the previous square determines a homotopy pushout in the Kan-Quillen model structure. This is true since the map $A_p\hookrightarrow B_p$ is a monomorphism.
\end{proof}

\begin{lemma}\label{preserved_pushout}
Consider a pushout square in $\sSet^{\text{fin}}$
\[
\begin{tikzcd}
A\arrow[d,hook]\arrow[r] & X\arrow[d] \\
B\arrow[r] & P\arrow[ul,"\ulcorner",very near start,phantom]
\end{tikzcd}
\]
where the map $A\hookrightarrow B$ is a monomorphism. Then the square
\[
\begin{tikzcd}
L(A)\arrow[d,hook]\arrow[r] & L(X)\arrow[d] \\
L(B)\arrow[r] & L(P)\arrow[ul,"\ulcorner",very near start,phantom]
\end{tikzcd}
\]
is a pushout in $\Cat$.
\end{lemma}

\begin{proof}
By definition, $L=(\text{Lan}_{\bar{y}} (\exit{\pos{-}}{\Delta^{-}}))\circ j$ where $j:\sSet^{\text{fin}}\hookrightarrow\Psh(\Delta)$ is fully faithful. By \Cref{pushout_sset}, $j$ sends the pushout square to a pushout square in $\Psh(\Delta)$. By \cite[Lemma 5.1.5.5]{htt}, left Kan extensions along Yoneda embeddings preserve colimits. In particular, it sends a pushout square in $\Psh(\Delta)$ to a pushout square in $\Cat$. Combining both observations shows that the square
\[
\begin{tikzcd}
L(A)\arrow[d,hook]\arrow[r] & L(X)\arrow[d] \\
L(B)\arrow[r] & L(P)\arrow[ul,"\ulcorner",very near start,phantom]
\end{tikzcd}
\]
is a pushout square in $\Cat$.
\end{proof}

The next step in proving \Cref{main_strat_B} is to show that $\exit{\pos{-}}{\norm{-}}$ also preserves some pushouts of finite simplicial sets. Namely, given a finite simplicial set $X$ and a top-dimensional nondegenerate simplex $\sigma$, we can decompose $X$ as a pushout 
\[
\begin{tikzcd}
\partial\Delta[\sigma]\arrow[r]\arrow[d,hook] & X'\arrow[d] \\
\Delta[\sigma]\arrow[r] & X\arrow[ul,"\ulcorner",phantom,very near start]
\end{tikzcd}.
\]

\begin{proposition}\label{pushout_exit}
Let $X$ be a finite simplicial set, $\sigma$ a top-dimensional nondegenerate simplex, and consider a pushout square

\[
\begin{tikzcd}
\partial\Delta[\sigma]\arrow[d,hook,"j"']\arrow[r] & X'\arrow[d,hook,"\iota"] \\
\Delta[\sigma]\arrow[r,"q"'] & X\arrow[ul,"\ulcorner",very near start,phantom]
\end{tikzcd}.
\]

Then the induced square
\[
\begin{tikzcd}
\exit{\pos{\partial\Delta[\sigma]}}{\partial\Delta^\sigma}\arrow[d]\arrow[r] & \exit{\pos{X'}}{\norm{X'}}\arrow[d] \\
\exit{\pos{\sigma}}{\Delta^\sigma}\arrow[r] & \exit{\pos{X}}{\norm{X}}
\end{tikzcd}
\]

is a pushout square of $\infty$-categories.
\end{proposition}

\begin{lemma}\label{pullback_sheaves}
Let $X$ be a finite simplicial set, $\sigma$ a top-dimensional nondegenerate simplex, and consider a pushout square

\[
\begin{tikzcd}
\partial\Delta[\sigma]\arrow[d,hook,"j"']\arrow[r] & X'\arrow[d,hook,"\iota"] \\
\Delta[\sigma]\arrow[r,"q"'] & X\arrow[ul,"\ulcorner",very near start,phantom]
\end{tikzcd}.
\]

Then the induced square
\[
\begin{tikzcd}
\Shv(\norm{X})\arrow[d]\arrow[r] & \Shv(\norm{X'})\arrow[d] \\
\Shv(\Delta^\sigma)\arrow[r] & \Shv(\partial\Delta^\sigma)
\end{tikzcd}
\]

is a pullback square of $\infty$-categories.
\end{lemma}

\begin{proof}

First observe that the square
\[
\begin{tikzcd}
\partial\Delta^{\sigma}\arrow[d,hook,"j"']\arrow[r] & \norm{X'}\arrow[d,"\iota",hook] \\
\Delta^{\sigma}\arrow[r,"q"'] & \norm{X}
\end{tikzcd}
\]
is a pullback square of topological spaces. This follows from the fact that both maps $j:\partial\Delta^{\sigma}\rightarrow\Delta^{\sigma},\iota:\norm{X'}\rightarrow \norm{X}$ are embeddings and an element $x\in\Delta^{\sigma}$ satisfies $q(x)\in\iota(\norm{X'})$ if and only if $x\in j(\partial\Delta^{\sigma})$. Moreover $q:\Delta^{\sigma}\rightarrow \norm{X}$ is proper since $\Delta^{\sigma}$ is compact and $\norm{X}$ is Hausdorff. This map also induces a homeomorphism $q:\Delta^{\sigma}\setminus\partial\Delta^{\sigma}=\overset{\circ}{\Delta}{}^{\sigma}\xrightarrow{\cong}\norm{X}\setminus \norm{X'}$. In particular, since $\norm{X'}$  is compact and $\norm{X}$ Hausdorff, this implies that the embedding $\iota:\norm{X'}\hookrightarrow \norm{X}$ is a closed embedding. Combining all these observations, we may apply \Cref{jansen_2_20} to get the desired result.
\end{proof}

\begin{lemma}\label{pullback_constructible_sheaves}
Let $X$ be a finite simplicial set, $\sigma$ a top-dimensional nondegenerate simplex, and consider a pushout square

\[
\begin{tikzcd}
\partial\Delta[\sigma]\arrow[d,hook,"j"']\arrow[r] & X'\arrow[d,hook,"\iota"] \\
\Delta[\sigma]\arrow[r,"q"'] & X\arrow[ul,"\ulcorner",very near start,phantom]
\end{tikzcd}.
\]

Then the induced square
\[
\begin{tikzcd}
\Shv^{\pos{X}}(\norm{X})\arrow[d]\arrow[r] & \Shv^{\pos{X'}}(\norm{X'})\arrow[d] \\
\Shv^{\pos{\sigma}}(\Delta^{\sigma})\arrow[r] & \Shv^{\pos{\partial\Delta[\sigma]}}(\partial\Delta^{\sigma})
\end{tikzcd}
\]

is a pullback square of $\infty$-categories.
\end{lemma}
\begin{proof}
By \Cref{pullback_sheaves}, it is enough to prove that any $F\in\Shv(\norm{X})$ pulling back to constructible sheaves is itself constructible. If $\tau$ is a stratum different from $\sigma$, then it is contained in $\norm{X'}$ and thus the restriction to this stratum is locally constant. For the $\sigma$-stratum, the restriction factors through the restriction of the sheaf pulled back to $\Delta^\sigma$ to the top stratum, which is locally constant. 
\end{proof}

\begin{proof}[Proof of \Cref{pushout_exit}]
This follows from \Cref{pullback_sheaves}, \Cref{pullback_constructible_sheaves} and \Cref{jansen_3_16}.
\end{proof}

\begin{proof}[Proof of \Cref{main_strat_B}]
It is enough to prove that the natural transformation $\eta:L\Rightarrow\exit{\pos{-}}{\norm{-}}$ provided by the universal property of the left Kan extension is a natural equivalence. 
Following the strategy of the proof of \cite[Proposition 2.2.2.7]{htt}, we do so by induction on the dimension of $X$ and on the number of top-dimensional nondegenerate simplices of $X$. First, we have $L(\emptyset)=\colim{\emptyset}{\exit{\pos{-}}{\Delta^-}}\simeq\emptyset\simeq\exit{\pos{\emptyset}}{\emptyset}$. Now let $X\in\sSet^{\text{fin}}$ and $\sigma$ a top-dimensional nondegenerate simplex in $X$. Then we can write $X$ as a pushout $X\cong X'\sqcup_{\partial\Delta[\sigma]}\Delta[\sigma]$. Now consider the cubical diagram in $\Cat$

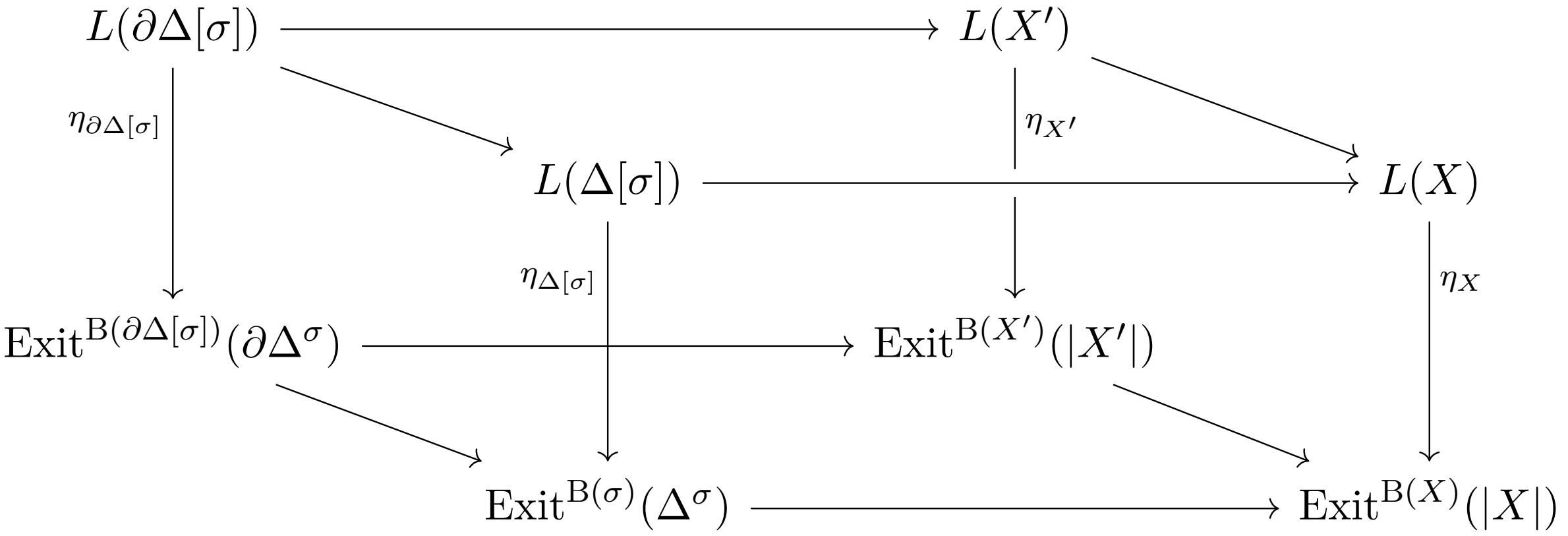
\[
\resizebox{\textwidth}{!}{%
\begin{tikzcd}[ampersand replacement=\&]
L(\partial\Delta[\sigma])
  \arrow[dd,"\eta_{\partial\Delta[\sigma]}"',near start]
  \arrow[rr]
  \arrow[dr]
\& \&
L(X')
  \arrow[dr]
  \arrow[dd,"\eta_{X'}",near start]
\&
\\
\&
L(\Delta[\sigma])
  \arrow[dd,"\eta_{\Delta[\sigma]}"',near start,crossing over]
  \arrow[rr,crossing over]
\& \&
L(X)
  \arrow[dd,"\eta_X",near start]
\\
\exit{\pos{\partial\Delta[\sigma]}}{\partial\Delta^{\sigma}}
  \arrow[rr]
  \arrow[dr]
\& \&
\exit{\pos{X'}}{\norm{X'}}
  \arrow[dr]
\&
\\
\&
\exit{\pos{\sigma}}{\Delta^\sigma}
  \arrow[rr]
\& \&
\exit{\pos{X}}{\norm{X}}
\end{tikzcd}%
}.
\]

By induction hypothesis, the functors $\eta_{\partial\Delta[\sigma]},\eta_{X'}$ are equivalences of $\infty$-categories. In addition, since $\Delta[\sigma]$ is representable, we have that $\eta_{\Delta[\sigma]}$ is also an equivalence. Moreover, by \Cref{preserved_pushout} and \Cref{pushout_exit}, the top and bottom faces are pushouts. These observations imply that the functor $\eta_{X}$ is also an equivalence, concluding the induction.

\end{proof}

We conclude this section by a conjecture. By \cite[Theorem A.6.10]{ha}, one gets for $p\geq 0$ an equivalence of $\infty$-categories
\[
\exit{\pos{p}}{\Delta^p}\xrightarrow{\simeq} N(\pos{p}).
\]
Moreover, the simplicial set $N(\pos{p})$ is isomorphic to the subdivision $\Sd(\Delta[p])$. Now the subdivision functor $\Sd:\sSet\rightarrow\sSet$ is defined by a left Kan extension. Thus, up to showing it is also a homotopy left Kan extension in the Joyal model structure, following the proof of \Cref{main_strat_B}, one would deduce a Joyal equivalence of simplicial sets
\[
\exit{\pos{X}}{\norm{X}}\simeq\Sd(X)
\]
for each finite simplicial set $X$. Combined with the exodromy equivalence, this motivates the following claim.

\begin{conjecture}
Let $X$ be a finite simplicial set. There is an equivalence of $\infty$-categories
\[
\Fun(\Sd(X),\cS)\simeq\Shv^{\pos{X}}{\norm{X}}.
\]
\end{conjecture}

\bibliographystyle{alpha}
\bibliography{bibliography_thesis}

\bigskip
\noindent
\textit{Email address:}
\href{mailto:rodrigue.haya@uclouvain.be}{\texttt{rodrigue.haya@uclouvain.be}}\\
\textsc{Institut de Recherche en Mathématique et Physique,
Université catholique de Louvain, chemin du cyclotron 2 bte L7.01.02, B–1348 Louvain-la-Neuve, Belgium}

\end{document}